\documentclass{amsart}
\usepackage[margin=1in]{geometry}
\usepackage{amsmath}
\usepackage{amssymb}
\usepackage{amsthm}
\usepackage{mathrsfs}
\usepackage{mathtools}
\usepackage{enumerate}
\usepackage{tikz}
\usetikzlibrary{decorations.pathreplacing}
\usepackage{hyperref}
\usepackage{tabularx}
\usepackage{float}
\usepackage[T1]{fontenc}
\newtheorem{thm}{Theorem}[section]
\newtheorem{prop}[thm]{Proposition}
\newtheorem{lem}[thm]{Lemma}
\newtheorem{cor}[thm]{Corollary}
\newtheorem*{thm*}{Theorem}
\theoremstyle{definition}
\newtheorem{dfn}[thm]{Definition}

\newtheorem{rmk}[thm]{Remark}
\newtheorem*{conj*}{Conjecture}
\newtheorem*{notat*}{Notation}

\numberwithin{equation}{section}

\usepackage{xcolor}

\DeclareMathOperator{\Vis}{Vis}
\DeclareMathOperator{\supp}{supp}
\DeclareMathOperator{\diam}{diam}
\DeclareMathOperator{\bbT}{\mathbb T}

\newcommand{\vs}{\vspace{.12in}}

\title{Improved Estimate for Visible Parts of Planar Ahlfors-regular Sets}

\author{Zhangze Li}
\address{Department of Mathematics, The Ohio State University, Columbus, Ohio 43210, USA}
\email{li.15128@osu.edu}

\author{Krystal Taylor}
\address{Department of Mathematics, The Ohio State University, Columbus, Ohio 43210, USA}
\email{taylor.2952@osu.edu}

\begin{document}

\begin{abstract}
    We study the visible parts of planar Ahlfors-regular sets. Orponen and Rutar recently proved that the visible part of a planar compact set has Hausdorff dimension at most $3/2$ in almost every direction, and that this bound is sharp. For compact $s$-Ahlfors-regular sets $E\subseteq\mathbb R^2$ with $1<s\leqslant 2$, we prove a refined estimate: the visible part of $E$ has Hausdorff dimension at most $\min\{s-\alpha_s(s-1),3/2\}$, where $\alpha_s=1/(3+\sqrt{8-2s})$, in almost every direction. This improves the previous bound for Ahlfors-regular sets due to D{\k{a}}browski. The main new ingredient is a quantitative Hausdorff content estimate for heavy slices of sub-uniformly distributed sets, obtained from the ball-tube incidence estimates of D{\k{a}}browski, Orponen, and Wang.
\end{abstract}

%s-(s-1)/(3+\sqrt{8-2s})

\maketitle

\setcounter{tocdepth}{1}
\tableofcontents

\section{Introduction and main result}
Let $\theta\in \mathbb S^{1}$. Consider the ``positive'' closed half line $\ell_{x,\theta}^+:=x+\mathbb R_{\geqslant 0}\theta$. 

\begin{dfn}[Visible Part]{}
    Let $E\subseteq\mathbb R^2$ be a compact set. The \emph{visible part of $E$ at direction $\theta$}, $\Vis_\theta(E)$, is the set of points $x\in E$ such that
    \[\ell_{x,\theta}^+\cap E=\{x\}.\]
\end{dfn}

For $\theta\in \mathbb S^1$, let $P_\theta:\mathbb R^2\to\mathbb R$ be the orthogonal projection on $\theta$, that is, $P_\theta(x)=x\cdot\theta$. Observe that $x\in\Vis_\theta(E)$ if and only if
\[y\cdot \theta\leqslant x\cdot \theta,\quad\forall\, y\in E\cap P_{\theta^{\perp}}^{-1}(P_{\theta^\perp}(x)).\]

We briefly describe the geometry of the visible part. Given $\theta\in \mathbb S^{1}$, the fibers of $P_{\theta^{\perp}}$ are lines parallel to $\ell_\theta$. These lines form a covering of $E$. On each line, we select the ``highest'' point in $E$, in the sense that there are no other points of $E$ beyond it in the direction $\theta$. $\Vis_\theta(E)$ is the collection of these ``highest'' points. It is called the ``visible part'' because if we look at $E$ from direction $\theta$, we can only see $\Vis_\theta(E)$.
\vs

By the Marstrand's projection theorem,
\[\dim_{\mathrm{H}}(P_{\theta^{\perp}}(E))=\min\{\dim_{\mathrm{H}} E,1\}\quad\text{for $\mathcal{H}^1$-a.e. $\theta\in \mathbb S^1$}.\]
Since $P_{\theta^{\perp}}(E)=P_{\theta^\perp}(\Vis_\theta(E))$, we have
\[\dim_{\mathrm{H}} \Vis_\theta(E)\geqslant \min\{\dim_{\mathrm{H}}E,1\}\quad\text{for $\mathcal{H}^1$-a.e. $\theta\in \mathbb S^{1}$}.\]

The visibility conjecture asserts that the equality should hold.
\begin{conj*}[Visibility Conjecture]{}
    Let $E\subseteq\mathbb R^2$ be a compact set. Then
    \[\dim_{\mathrm{H}} \Vis_\theta(E)=\min\{\dim_{\mathrm{H}}E,1\}\quad\text{for $\mathcal{H}^{1}$-a.e. $\theta\in \mathbb S^{1}$}.\]
\end{conj*}

There are many partial results towards this conjecture in the past few years. It is confirmed for fractal percolation \cite{AJJRS12}, several classes of self-similar sets \cite{FF13,JJMON03,JJN04,JJSW22,KMP26}, and
dominated self-affine sets \cite{ER21}. For general compact sets, the upper bound $\dim_{\mathrm{H}} \Vis_\theta(E)\leqslant 1.99$ was proved by Orponen \cite{Tuo23}, and later D{\k{a}}browski \cite{Dabrowski-vis24} improved it to $\dim_{\mathrm{H}} \Vis_\theta(E)\leqslant 11/6$. The visibility conjecture can also be formulated and studied in higher dimensions; see \cite{Dabrowski-vis24,HJL-higherdim,Tuo23}.
\vs

However, a recent breakthrough from Orponen and Rutar \cite{TR26} disproved this conjecture. They showed that the Hausdorff dimension of visible part is at most $3/2$ in almost every direction, and this is sharp.
\begin{thm*}[\cite{TR26}]
    Let $E\subseteq\mathbb R^2$ be a compact set. Then 
    \[\dim_{\mathrm H}\Vis_\theta(E)\leqslant \frac{3}{2} \quad\text{for $\mathcal{H}^1$-a.e.}\ \theta\in\mathbb S^1.\]
    In particular, there exists a compact set such that its visible part has Hausdorff dimension $3/2$ for every $\theta\in\mathbb S^1$. 
\end{thm*}

However, the sharp example constructed in \cite{TR26} is neither Ahlfors-regular nor connected. One may ask whether the visibility conjecture holds for Ahlfors-regular or connected sets? O'Neil \cite{ONeilconnected} established a very strong bound for connected compact sets. Meanwhile, D{\k{a}}browski \cite{Dabrowski-vis24} has proved that if $E\subseteq\mathbb R^2$ is a compact $s$-Ahlfors-regular set with $1<s\leqslant 2$, then
\begin{equation}\label{eqn:Dabrowski bound}
    \dim_{\mathrm H}\Vis_\theta(E)\leqslant s-\alpha(s-1) \quad\text{for $\mathcal{H}^1$-a.e.}\ \theta\in\mathbb S^1,
\end{equation}
where $\alpha=1/(3+\sqrt{6})$. Note that when $s=3/2$, $\dim_{\mathrm H}\Vis_\theta(E)<3/2$ for $\mathcal{H}^1$-a.e.\ $\theta\in\mathbb S^1$.
\vs

Our main result improves \eqref{eqn:Dabrowski bound}. We first recall the definition of Ahlfors-regularity.
\begin{dfn}[Ahlfors-regular set]{}
    Let $E\subseteq \mathbb R^2$ be a compact set and $0<s\leqslant 2$. We say that $E$ is \emph{$s$-Ahlfors-regular} if it supports a finite nonzero Borel measure $\mu$ satisfying
    \[C^{-1}r^s\leqslant \mu(B(x,r))\leqslant Cr^s \quad\forall\,x\in E,\ 0<r\leqslant\diam E.\]
\end{dfn}

\begin{thm}[Main theorem]\label{thm:main result}
    Let $E\subseteq\mathbb R^2$ be a compact $s$-Ahlfors-regular set with $1<s\leqslant 2$. Then
    \begin{equation}\label{eqn:main bound}
        \dim_{\mathrm H}\Vis_\theta(E)\leqslant \min\left\{s-\alpha_s(s-1),\frac{3}{2}\right\} \quad\text{for $\mathcal{H}^1$-a.e.}\ \theta\in\mathbb S^1,
    \end{equation}
    where $\alpha_s=1/(3+\sqrt{8-2s})$.
\end{thm}

%\begin{thm}
    %The conclusion holds if its measure hypothesis is replaced by the following one: there is a Borel probability measure $\mu$ with support $E$ such that, for every $u>0$, some $C_u\geqslant 1$ satisfies
    %\[C_u^{-1}r^{s+u}\leqslant \mu(B(x,r))\leqslant C_u r^{s-u}\quad\forall\,x\in E,\ 0<r\leqslant 1. \]
    %The constants are uniform in $x$ and $r$. In particular,
    %$\dim E=\overline{\dim}_{B}E=s$.
%\end{thm}

We briefly discuss the strategy to prove Theorem \ref{thm:main result}. We follow the visibility decomposition developed in \cite{Dabrowski-vis24,Tuo23}. 
After normalizing, let $E\subseteq\mathbb [0,1)^2$ be a compact $s$-Ahlfors-regular set. 
Our aim is to obtain, for a suitable $\tau>0$,
\[\mathcal{H}_{\infty}^{s-\tau}( \Vis_\theta(E) )=0\quad\text{for $\mathcal{H}^1$-a.e.\ $\theta\in\mathbb S^1$}.\]

We first fix two dyadic scales, a fine scale $\delta$ and a coarse scale $\Delta\gg \delta$, and we consider the corresponding grids of dyadic cubes in $[0,1)^2$. 
For a fixed direction $\theta$, we cover the unit square by parallel tubes $T$ of width $\delta$. Within each tube, we count the fine $\delta$-cubes meeting $E\cap T$, grouping them according to the coarse $\Delta$-cubes that contain them.   
\vs

After estimating a light part, we separate the remaining points into good, bad, and heavy parts, estimating each contribution separately.
For the good-part estimate, our goal is to use fewer fine cubes than are needed to cover all of $E$. By identifying blocked or hidden points, we can omit them from the cover. 
\vs

If a tube $T$ meets the good part, then it is neither light nor contained in a globally heavy tube. Thus, it contains many occupied $\delta$-cubes, distributed among a controlled number of coarse $\Delta$-cubes. 
By the pigeonhole principle, there exists a $\Delta$-cube $Q$ containing a sufficiently large number of these $\delta$-cubes. We say that $T$ is rich relative to that $Q$, see Lemma \ref{lem:rich tube} for the details.  
\vs

%For each good point in a tube $T$, its viewing line meets $E\cap Q$ if $T$ is relative rich to $Q$. Thus, good points lying behind this intersection are not visible, and can be omitted from the cover. 
Among the cubes relative to which $T$ is rich, choose $Q_T$ to maximize $\inf\{x\cdot\theta:x\in E\cap3Q_T\}$. Counting the good points that may remain visible then gives the estimate in Proposition \ref{prop:delta covering of good part}.
\vs

The difficulty is that richness of a tube does not guarantee that each viewing line inside it intersects $E\cap Q$. A line can pass through a tube rich relative to $Q$ but miss $E\cap Q$, leaving points behind the cube potentially visible. These are the bad lines.  
A Fourier-analytic argument gives a Hausdorff content upper bound for projections of bad lines in terms of the Sobolev norms of the projected measures; see Proposition \ref{prop:bad lines estimate}. 
Outside the global heavy part, the projection covers can be lift to cover the bad part of $E$. Combining this with the averaged Sobolev estimates from Marstrand's projection theorem gives the bad-part estimate in Corollary \ref{cor:bad part}.
\vs

Finally, we handle the heavy part, consisting of points of $E$ lying in tubes whose populations exceed specified covering-number thresholds. Global heaviness measures the population along a tube at the scale of its width, while local heaviness measures the fine $\delta$-cube population inside a coarse $\Delta$-cube. These conditions model, at finite scales, the exceptional slicing behavior suggested by Marstrand's slicing theorem. 
This theorem states that if $A\subseteq\mathbb R^2$ is a compact set with $0<\mathcal{H}^t(A)<\infty$ for $1<t\leqslant 2$, then 
\[\dim_{\mathrm{H}} (A\cap \ell_{x,\omega})=t-1\quad\text{for $\mathcal{H}^t\times\mathcal{H}^1$-a.e.}\ (x,\omega)\in A\times\mathbb S^1,\]
where the line $\ell_{x,\omega}:=x+\mathbb R\omega$. 
Thus, for almost every direction, the points lying on slices of dimension greater than $t-1$ form a set of zero $\mathcal H^t$-measure. 
\vs 

A maximal function estimate established by D{\k{a}}browski \cite{Dabrowski-vis24} combined with the averaged Sobolev estimates handles the global heavy part in Corollary \ref{cor:global heavy part}. 
\vs

Our improvement concerns the local heavy part. Instead of applying the maximal-function estimate used in \cite{Dabrowski-vis24}, we prove a quantitative Hausdorff-content estimate using the ball-tube incidence estimate of D{\k{a}}browski, Orponen, and Wang \cite{DOW24-heavyslice}.
Applying this estimate to rescaled coarse cubes gives the local heavy-part estimate in Proposition \ref{prop:local heavy part}, which leads to the bound in \eqref{eqn:main bound}.
\vs

The organization is as follows. 
Section \ref{sec:prelim} contains preliminaries on Frostman's lemma, Sobolev norm estimates, and maximal function estimates. Section \ref{sec:incidence} discusses Marstrand's slicing theorem and proves a Hausdorff content estimate on heavy slices. Section \ref{sec:proof of main} is the proof of Theorem \ref{thm:main result}.
\vs

\noindent {\bf Acknowledgment}.
K.T.\ is supported in part by the Simons Foundation Grant GR137264. The authors wish to thank Damian D{\k{a}}browski for many helpful comments and discussions on this article.

\medskip
\noindent {\bf AI declaration}.
All ideas and arguments presented in this manuscript are original contributions of the authors. AI assistance was used only to identify typographical errors. The manuscript was written entirely by the authors.

\section{Preliminaries}\label{sec:prelim}
\subsection{Frostman's lemma and Sobolev norm estimate}
One useful tool to calculate the Hausdorff dimension is the following Frostman's lemma. The proof can be found in \cite[Theorem~2.7]{Mattilabook}.

\begin{lem}[Frostman's lemma]{}
    Let $0<s\leqslant 2$, and $E\subseteq\mathbb R^2$ be a Borel set with $\mathcal{H}_\infty^s(E)>0$. Then there exists a Borel probability measure $\mu$ supported on $E$ such that
    \[\mu(B(x,r))\lesssim \mathcal{H}_\infty^s(E)^{-1} r^s,\quad \forall\,x\in\mathbb R^2,r>0.\]
\end{lem}

%The second variant Frostman's lemma is from \cite{Tuo23}, which provides a lower bound on the measure on the dyadic cubes.

%\begin{lem}[Lemma~A.1 in \cite{Tuo23}]\label{lem:Orponen's Frostman lemma}
    %Let $0\leqslant s\leqslant 2$, and $E\subseteq[0,1)^2$ be a compact set. Then there exists a Radon measure $\mu$ supported on $E$ such that
    %\[\mu(B(x,r))\lesssim r^s,\quad \forall\,x\in\mathbb R^2,\ r>0,\]
    %and
    %\[\mu(\overline{Q})\gtrsim \min\{\mathcal{H}^{s}_\infty(E\cap Q),|Q|\},\quad\forall\,Q\in\mathcal{D}([0,1)^2).\]
%\end{lem}

We shall need the Sobolev norm when dealing with the orthogonal projection. Let $\mu$ be a finite Borel measure in $\mathbb R^2$. For $\theta\in \mathbb S^1$, define the pushforward measure $\mu_\theta:={{P_\theta}}_\sharp\mu$.

\begin{dfn}[Sobolev norm]
    Let $s\in\mathbb R$, and $f:\mathbb R^n\to \mathbb C$ be a Schwartz function. Then the \emph{$s$-homogeneous Sobolev norm of $f$} is defined as
    \[\|f\|_{\dot H^s(\mathbb R^n)}:=\left(\int_{\mathbb R^n} \left|\widehat{f}(\xi)\right|^2|\xi|^{2s}\,d\xi\right)^{1/2}.\]
    The \emph{$s$-inhomogeneous Sobolev norm of $f$} is defined as
    \[\|f\|_{H^s(\mathbb R^n)}:=\left(\int_{\mathbb R^n} \left|\widehat{f}(\xi)\right|^2\left(1+|\xi|^2\right)^{s}\,d\xi\right)^{1/2}.\]
\end{dfn}

Recall that for $0<s<n$ and a finite Borel measure $\mu$ on $\mathbb R^n$, the \emph{$s$-Riesz energy of $\mu$} is defined as 
\[I_s(\mu):=\int_{\mathbb R^n} |x-y|^{-s}\,d\mu(x)\,d\mu(y).\]
It is well-known that 
\[I_s(\mu)\sim_{s,n} \|\mu\|_{\dot H^{(s-n)/2}(\mathbb R^n)}^2,\]
which allows the connection between the Riesz energy and the Sobolev norm. This leads to the following proposition.

\begin{prop}[Sobolev smoothing under projections]{}\label{prop:sobolev proj}
    Let $\mu$ be a finite Borel measure supported on $[0,1]^2$ and $1<s\leqslant 2$. Suppose $\mu(B(x,r))\leqslant Cr^s$ for every $x\in\supp \mu$ and $r>0$. Then for every $0<\varepsilon<s-1$, we have
    \[\int_{\mathbb S^1} \|\mu_{\theta}\|_{H^{(s-1-\varepsilon)/2}(\mathbb R)}^2\,d\mathcal{H}^1(\theta)\lesssim_\varepsilon C\mu(\mathbb R^2).\]
\end{prop}

Proposition \ref{prop:sobolev proj} can be used to prove the Marstrand's projection theorem. We sketch the proof here. See \cite[Theorem~5.10]{Mattilabook} for a detailed proof. 

\begin{proof}
    We first note that for every $t\geqslant 0$, we have
    \begin{align*}
        \|\mu_{\theta}\|_{H^t(\mathbb R)}^2&=\int_{\mathbb R} |\widehat{\mu_{\theta}}(\xi)|^2(1+|\xi|^2)^t\,d\xi\\
        &\lesssim \int_{|\xi|<1} |\widehat{\mu_{\theta}}(\xi)|^2\,d\xi+\int_{|\xi|\geqslant 1} |\widehat{\mu_{\theta}}(\xi)|^2 |\xi|^{2t}\,d\xi\\
        &\leqslant \mu(\mathbb R^2)^2+\|\mu_{\theta}\|_{\dot H^t(\mathbb R)}^2\\
        &\lesssim C\mu(\mathbb R^2)+\|\mu_{\theta}\|_{\dot H^t(\mathbb R)}^2.
    \end{align*}
    The last step uses $\mu(\mathbb R^2)=\mu([0,1]^2)\leqslant\mu(B(0,\sqrt{2}))\lesssim C$.

    Thus, it suffices to estimate $\|\mu_\theta\|_{\dot H^{(s-1-\varepsilon)/2}(\mathbb R)}^2$. Since $\widehat{\mu_{\theta}}(\xi)=\widehat{\mu}(\xi\theta)$ for every $\xi\in\mathbb R$ and $\widehat{\mu}(y)=\overline{\widehat{\mu}(-y)}$ for every $y\in\mathbb R^2$, we have
    \begin{align*}
        \int_{\mathbb S^1} \|\mu_{\theta}\|_{\dot H^{(s-1-\varepsilon)/2}(\mathbb R)}^2\,d\mathcal{H}^1(\theta)&=2\int_{\mathbb S^1} \int_0^\infty|\widehat{\mu}(\xi\theta)|^2 \xi^{s-1-\varepsilon}\,d\xi\,d\mathcal{H}^1(\theta)\\
        &\sim \int_{\mathbb R^2} |\widehat{\mu}(y)|^2|y|^{s-2-\varepsilon}\,dy\\
        &\sim_\varepsilon I_{s-\varepsilon}(\mu)\\
        &\lesssim_\varepsilon C\mu(\mathbb R^2).
    \end{align*}
\end{proof}

\subsection{Maximal function estimate}
Let $\nu$ be a finite Borel measure on $\mathbb R^n$. The \emph{Hardy-Littlewood maximal function of
$\nu$} is defined as
\[\mathcal{M}\nu(x):=\sup_{r>0} r^{-n}\nu(B(x,r)).\]

The following maximal function estimate is established in \cite{Dabrowski-vis24}. It nicely relates the heavy slices with the Sobolev norm of push-forward of the measure under the orthogonal projection. It can be viewed as a quantitative version of the exceptional set estimate in Marstrand's slicing theorem. We discuss this exceptional set estimate in more detail in Section \ref{sec:incidence}.

\begin{prop}[Proposition~3.1 in \cite{Dabrowski-vis24}]{}\label{prop:maximal function estimate}
    Let $1<s\leqslant 2$ and $0<\varepsilon<(s-1)/2$. Let $\mu$ be an $s$-Ahlfors-regular measure supported in $B(0,1)$. For every $M>10\mu(\mathbb R^2)$, let
    \[H_{\theta,M}:=\{x\in\supp\mu:\mathcal{M}\mu_{\theta}(P_{\theta}(x))\geqslant M\}.\]
    Then
    \[\mathcal{H}^{1+2\varepsilon}_\infty(H_{\theta,M})\lesssim_\varepsilon M^{-1}\|\mu_{\theta}\|_{H^{(s-1-\varepsilon)/2}(\mathbb R)}^2.\]
\end{prop}

\section{Heavy slice estimate for sub-uniformly distributed sets}\label{sec:incidence}
Recall the Marstrand's slicing theorem: if $E\subseteq\mathbb R^2$ is a compact set with $0<\mathcal{H}^s(E)<\infty$ for $1<s\leqslant 2$, then 
\[\dim_{\mathrm{H}} (E\cap \ell_{x,\theta})=s-1\quad\text{for $\mathcal{H}^s\times\mathcal{H}^1$-a.e.}\ (x,\theta)\in E\times\mathbb S^1.\]
Mattila \cite[Section~6.4]{Mattilabook} asked whether one can obtain an upper bound for the Hausdorff dimension of the exceptional set of heavy slices
\[\{x\in E:\dim_\mathrm{H}(E\cap \ell_{x,\theta})\geqslant t\},\quad s-1<t\leqslant 1.\]

This problem was recently resolved using Furstenberg set estimates; see \cite[Theorem~1.7 \& Theorem~1.8]{DOW24-heavyslice}.
\begin{thm}[\cite{DOW24-heavyslice}]{}
    Let $E\subseteq\mathbb R^2$ be a compact set with $\dim_\mathrm{H}E=s\in [1,2]$. The following estimates are sharp.
    
    If $\frac{1}{3}(2s-1)<t\leqslant 1$, then $\{x\in E:\dim_\mathrm{H}(E\cap \ell_{x,\theta})\geqslant t\}=\varnothing$ for $\mathcal{H}^1$-a.e.\ $\theta\in\mathbb S^1$. If $s-1<t \leqslant\frac{1}{3}(2s-1)$, then $\dim_\mathrm{H}\{x\in E:\dim_\mathrm{H}(E\cap \ell_{x,\theta})\geqslant t\}\leqslant \min\{2s-3t,s\}$ for $\mathcal{H}^1$-a.e.\ $\theta\in\mathbb S^1$. In particular, for $\mathcal{H}^1$-a.e.\ $\theta\in\mathbb S^1$, 
    \[\dim_\mathrm{H}\{x\in E:\dim_\mathrm{H}(E\cap \ell_{x,\theta})>s-1\}\leqslant \min\{3-s,s\}.\]
\end{thm}

For Ahlfors-regular sets, D{\k{a}}browski \cite{Dabrowski-vis24} obtained a stronger estimate using Proposition \ref{prop:maximal function estimate}.
\begin{thm}[\cite{Dabrowski-vis24}]
    Let $E\subseteq\mathbb R^2$ be a compact $s$-Ahlfors-regular set with $1<s\leqslant 2$. Then for $\mathcal{H}^1$-a.e.\ $\theta\in\mathbb S^1$, 
    \[\dim_\mathrm{H}\{x\in E:\overline{\dim}_\mathrm{M}(E\cap \ell_{x,\theta})>s-1\}\leqslant 1.\]
\end{thm}

D{\k{a}}browski, Orponen, and Wang \cite{DOW24-heavyslice} subsequently extended the above result to a broader class called sub-uniformly distributed sets.
\begin{thm}[\cite{DOW24-heavyslice}]{}\label{thm:DOW heavy slice bound}
    Let $E\subseteq\mathbb R^2$ be a compact sub-uniformly distributed set with $\dim_\mathrm{H}E=s\in [1,2)$. Then for $s-1<t\leqslant 1$, 
    \[\dim_\mathrm{H}\{x\in E:\dim_\mathrm{H}(E\cap \ell_{x,\theta})\geqslant t\}\leqslant s-t \quad\text{for $\mathcal{H}^1$-a.e.\ $\theta\in\mathbb S^1$}.\]
\end{thm}

We will use a ball-tube incidence version of Theorem \ref{thm:DOW heavy slice bound} in the proof of Theorem \ref{thm:main result}. The purpose of this section is to establish a quantitative Hausdorff content estimate for heavy slices; see Proposition \ref{prop:local heavy from incidence}. 

\subsection{Sub-uniformly distributed sets}
We first introduce the definition and basic properties of sub-uniformly distributed sets.

\begin{notat*}[Dyadic cubes and covering number]{}
    Let $\delta\in 2^{-\mathbb N}$. We use $\mathcal{D}_\delta(\mathbb R^2)$ to represent all dyadic $\delta$-cubes in $\mathbb R^2$. Let $E\subseteq\mathbb R^2$ be a bounded set. Then define
    \[\mathcal{D}_\delta(E)=E_\delta:=\{Q\in \mathcal{D}_\delta(\mathbb R^2):E\cap Q\ne\varnothing\}.\]
    We abbreviate $\mathcal{D}_\delta:=\mathcal{D}_\delta([0,1)^2)$.
    We also use $|E|_\delta:=|\mathcal{D}_\delta(E)|$ to represent dyadic $\delta$-covering number of $E$.
\end{notat*}

\begin{dfn}[Sub-uniformly distributed set]
    Let $E\subseteq\mathbb R^2$ be a bounded set. We say that $E$ is \emph{sub-uniformly distributed} if there exists a constant $C>0$ such that
    \[|E|_R\cdot |E\cap Q|_r\leqslant C|E|_r,\quad\forall\,Q\in\mathcal{D}_R(E),\ 0<r\leqslant R<\infty.\]
\end{dfn}

It is clear that if $E$ is $s$-Ahlfors-regular, then $E$ is sub-uniformly distributed, because $|E\cap Q|_r\lesssim (R/r)^s\sim |E|_r/|E|_R$ for every $Q\in\mathcal{D}_R(E)$ and $0<r<R<\diam E$.

A useful observation is that, for sub-uniformly distributed sets, the lower Minkowski dimension coincides with the Hausdorff dimension.
\begin{lem}[Lower Minkowski dimension]{}\label{lem:subuniform Mink equal to Haus}
    Let $E\subseteq\mathbb R^2$ be a compact sub-uniformly distributed set. Then $\dim_{\mathrm{H}} E=\underline{\dim}_\mathrm{M} E$.
\end{lem}

We now give a $\delta$-discretized definition for sub-uniformly distributed sets.
%We need the following $\delta$-discretized definition for sub-uniformly distributed sets, which allows us to decompose a collection of $\delta$-dyadic cubes into sub-uniformly distributed subcollections.
\begin{dfn}[$\delta$-discretized definition]
     Let $C>0$ and $\{\Delta_n<\Delta_{n-1}<\dots<\Delta_0\}\subseteq 2^{-\mathbb N}$ be an increasing sequence of scales. Let $\mathcal{P}$ be either a subset of $B(0,1)$ or a subfamily of $\mathcal{D}_\delta$ with $\delta\leqslant\Delta_n$. We say $\mathcal{P}$ is \emph{sub-$\{\Delta_j\}_{j=0}^n$-uniform with constant $C$} if
     \[|\mathcal{P}|_{\Delta_j}\cdot|\mathcal{P}\cap p|_{\Delta_{j+1}}\leqslant C|\mathcal{P}|_{\Delta_{j+1}},\quad\forall\, p\in\mathcal{D}_{\Delta_j}(\mathcal{P}),\ j\in\{0,\dots,n-1\}.\]
\end{dfn}

With the $\delta$-discretized definition, the following decomposition lemma allows us to decompose a collection of $\delta$-dyadic cubes into sub-uniformly distributed subcollections. We will use this decomposition later in the incidence estimate for tubes.
\begin{lem}[Decomposition]{}\label{lem:decomp}
    Let $0<\eta<1$ with $\eta^{-1}\in\mathbb N$, and set $S:=\delta^{-\eta}\in 2^\mathbb N$. Consider the scale sequence 
    $\{\Delta_j=S^{n-j}\delta:j\in\{0,1,\dots,n\}\}\subseteq 2^{-\mathbb N}$, where $n=\eta^{-1}$. Then for every $\mathcal{P}\subseteq\mathcal{D}_\delta$, there exists a partition $\mathcal{P}=\bigsqcup\limits_{i=1}^N\mathcal{P}_i$ such that $\mathcal{P}_i$ are pairwise disjoint,  each $\mathcal{P}_i$ is sub-$\{\Delta_j\}_{j=0}^n$-uniform with constant $2$, and $N\leqslant (-Cn^{-1}\log \delta)^n$.
\end{lem}

The proofs of Lemma \ref{lem:subuniform Mink equal to Haus} and \ref{lem:decomp} can be found in \cite{DOW24-heavyslice}.

\subsection{Incidence estimate}
In this section, we use ball-tube incidence estimates to  establish a quantitative Hausdorff content bound for heavy slices. We begin by introducing some notation from incidence geometry.
\subsubsection{Incidence geometry notation}
\begin{notat*}[Dyadic tubes]{}
    Let $\delta\in 2^{-\mathbb N}$. Consider the line
    \[\mathbf{D}(a,b):=\{(x,y)\in\mathbb R^2:y=ax+b\}\]
    The collection of dyadic $\delta$-tubes is defined as
    \[\mathcal{T}^\delta:=\{T:\exists\,p\in\mathcal{D}_\delta,\ T =\mathbf D(p)\}.\]
\end{notat*}

\begin{notat*}[Incidence]
    Let $\delta\in 2^{-\mathbb N}$. For every $p\in\mathcal{D}_\delta$, define the associated ball $B_p$ that has the same center as $p$ with $\diam B_p=10\delta$. Similarly, for every $T=\mathbf{D}([a,a+\delta)\times[b,b+\delta))\in\mathcal{T}^\delta$, define the associated standard tube $\bbT_T$ that has the center line $\mathbf{D}(a,b)$ with width $10\delta$. 
    
    Let $\mathcal{P}\subseteq \mathcal{D}_\delta$ and $\mathcal{T}\subseteq\mathcal{T}^\delta$. The incidence between $\mathcal{P}$ and $\mathcal{T}$ is defined as
    \[\mathcal{I}(\mathcal{P},\mathcal{T}):=\{(p,T):B_p\cap\bbT_T\ne\varnothing\}.\]

    Note that we have $T\cap B(0,1)\subseteq\bbT_T$. Also, if $p,\mathbf p\in \bigcup\limits_{\delta}\mathcal{D}_\delta$ with $p\subseteq \mathbf p$ and $T,\mathbf T\in \bigcup\limits_{\delta}\mathcal{T}^\delta$ with $T\subseteq\mathbf{T}$, then $B_p\subseteq B_{\mathbf p}$ and $\bbT_T\cap B(0,2)\subseteq \bbT_{\mathbf{T}}\cap B(0,2)$. Thus, if $B_p\cap \bbT_T\ne \varnothing$, then $B_{\mathbf{p}}\cap \bbT_{\mathbf{T}}\ne\varnothing$.
\end{notat*}

The dyadic tubes are convenient for the dyadic decomposition and calculation, while the standard tubes depict the geometry of the incidence.

The following containment lemma states that if $(p,T)$ is viewed as dyadic cubes in $ \mathcal{D}_\delta([0,1)^4)$, then we have a containment of incidence for coarser scales.
\begin{lem}[Scale containment]\label{lem:containment}
    Let $\delta\in 2^{-\mathbb N}\cap (0,1/100]$, $\mathcal{P}\subseteq \mathcal{D}_\delta$ and $\mathcal{T}\subseteq \mathcal{T}^\delta$. Then
    \[(\mathcal{I}(\mathcal{P},\mathcal{T}))_\Delta\subseteq \mathcal{I}(\mathcal{P}_\Delta,\mathcal{T}_\Delta),\quad\forall\,\Delta\in 2^{-\mathbb N}\cap [\delta,1]\]
\end{lem}

\subsubsection{Hausdorff content estimate}
The main tool we use in the proof of Theorem \ref{thm:main result} is the following incidence estimate, which serves as a $\delta$-discretized version of Theorem \ref{thm:DOW heavy slice bound}.
\begin{prop}[Proposition~3.2 in \cite{DOW24-heavyslice}]{}\label{prop:DOW incidence}
    Let $C>0$, $0<\kappa<1$ and $0<\eta\leqslant\kappa/2$. The following holds for all sufficiently small $\delta\in 2^{-\mathbb N}$. Set $S:=\delta^{-\eta}\in 2^{\mathbb N}$. Consider the scale sequence 
    $\{\Delta_j=S^{n-j}\delta:j\in\{0,1,\dots,n\}\}$. Let $\mathcal{T}\subseteq\mathcal{T}^\delta$ and $\mathcal{P}\subseteq\mathcal{D}_\delta$ be sub-$\{\Delta_j\}_{j=0}^n$-uniform with constant $C$. Suppose that
    \[|\mathcal{I}(\mathcal{P},\mathcal{T})|\geqslant \delta^{1-\kappa}|\mathcal{P}||\mathcal{T}|.\]
    Then there exists a scale $\Delta\in\{\Delta_n,\dots,\Delta_1\}\subseteq [\delta,\delta^\eta]$ such that
    \[|\mathcal{I}(\mathcal{P}_\Delta,\mathcal{T}_\Delta)|\lesssim_{C,\eta}\delta^{\kappa-2\eta}\Delta^{-2},\]
    and
    \[\frac{|\mathcal{I}(\mathcal{P}_\Delta,\mathcal{T}_\Delta)|}{\Delta|\mathcal{P}_\Delta||\mathcal{T}_\Delta|}\gtrsim_{C,\eta}\delta^{\eta-\kappa}.\]
    %Here $\mathcal{P}_\Delta:=\mathcal{D}_\Delta(\mathcal{P})$ and $\mathcal{T}^\Delta(\mathcal{T})$ are the dyadic $\Delta$-covers of $\mathcal{P}$ and $\mathcal{T}$, respectively.
\end{prop}

The proof of Proposition \ref{prop:DOW incidence} is based on the ``high-low method''. Roughly speaking, we cut the incidence into ``high frequency term'' and ``low frequency term''. Sub-uniform distribution assumption guarantees that ``high frequency term'' will eventually occur at a scale. See \cite[Proposition~3.2]{DOW24-heavyslice} for a detailed proof.

Before stating the quantitative estimate on heavy slices, we introduce the following notation.
\begin{notat*}[Single tube incidence and direction set]
    Given $\mathcal{P}\subseteq\mathcal{D}_\delta$ and $T=\mathbf{D}([a,a+\delta)\times [b,b+\delta))\in\mathcal{T}^\delta$, define the incidence set of $T$ as
    \[\mathcal{I}(\mathcal{P},\{T\}):=\{p\in \mathcal{P}:B_p\cap \bbT_T\ne\varnothing\},\]
    and the direction set of $T$ as
    \[\Theta_T:=\left\{\pm \frac{(1,u)}{\sqrt{1+u^2}}:u\in [a,a+\delta)\right\}\subseteq\mathbb S^1.\]
\end{notat*}

Our main estimate that will be used later in Section \ref{sec:proof of main} is the following. It converts a lower bound on the single tube incidence into an averaged Hausdorff content estimate on the heavy slices.
Informally, it says that 
requiring many incidences per tube forces the covered portion of $E$ to have small Hausdorff content on average over directions.
\begin{prop}[From incidence to averaged Hausdorff content]{}\label{prop:local heavy from incidence}
    Let $E\subseteq[0,1)^2$ be a compact sub-uniformly distributed set with $\dim_{\mathrm H} E=s\in (1,2)$. Let $1<t\leqslant s,0<\kappa<1$ and $\varepsilon>0$.

    The following holds for all sufficiently small $\delta\in 2^{-\mathbb N}$. Let $\mathcal{P}:=\mathcal{D}_\delta(E)$, $\mathcal{T}\subseteq \mathcal{T}^\delta$.
    %For each $T=\mathbf D([a,a+\delta)\times[b,b+\delta))\in\mathcal{T}$, define the incidence set of $T$ as
    %\[\mathcal{I}(\mathcal{P},\{T\}):=\{p\in \mathcal{P}:B_p\cap \bbT_T\ne\varnothing\},\]
    %and the direction set of $T$ as
    %\[\Theta_T:=\left\{\pm \frac{(1,u)}{\sqrt{1+u^2}}:u\in [a,a+\delta)\right\}\subseteq\mathbb S^1.\]
    %Also, define heavy part of $E$ at direction $\theta$ as
    %\[E_{\theta, H}:=E\cap \bigcup_{T\in\mathcal{T}:\theta\in\Theta_{T}}\bbT_T.\]
    Suppose that
    \begin{equation}\label{eqn:heavy incidence assumption}
        |\mathcal{I}(\mathcal{P},\{T\})|\geqslant\delta^{1-\kappa}|\mathcal{P}|, \quad\forall\,T\in\mathcal{T}.
    \end{equation}
    For $\theta\in\mathbb S^1$, define the $\theta$-directional tube component of $E$ as
    \[E_{\theta,\mathcal{T}}:=E\cap \bigcup_{T\in\mathcal{T}:\theta\in\Theta_{T}}\bbT_T.\]
    Then
    \[\int_{\mathbb S^1} \mathcal{H}^{t}_\infty(E_{\theta, \mathcal{T}})\,d\mathcal{H}^1(\theta)\lesssim_{s,t,\kappa,\varepsilon} \delta^{\frac{\kappa(t+1-s)}{2-s}-\varepsilon}.\]
\end{prop}

\begin{proof}
    By Lemma \ref{lem:subuniform Mink equal to Haus}, $\dim_{\mathrm H} E=\underline{\dim}_{\mathrm M} E=s$. Thus, for every $v>0$, there exists a $C_v>0$ such that $|E|_\Delta\geqslant C_v \Delta^{-s+v}$ for every $0<\Delta\leqslant 1$.

    Since $t>1$ and $s<2$, we may choose $v=v(s,t,\kappa,\varepsilon)>0$ sufficiently small so that
    \[\kappa(t-1)\left(\frac{1}{2-s}-\frac{1}{2-s+v}\right)<\varepsilon/3.\]
    Next, we choose $\eta=\eta(s,t,\kappa,\varepsilon)\in (0,\kappa/2)$ sufficiently small so that
    \[\eta\left(2+\frac{t-1}{2-s+v}\right)<\varepsilon/3.\]

    Set $S:=\delta^{-\eta}$. Consider the scale sequence $\{\Delta_j=S^{n-j}\delta:j\in\{0,1,\dots,n\}\}$. By Lemma \ref{lem:decomp}, there exists a partition $\mathcal{T}=\bigsqcup\limits_{i=1}^N \mathcal{T}_i$ such that $\mathcal{T}_i$ are pairwise disjoint, each $\mathcal{T}_i$ is sub-$\{\Delta_j\}_{j=0}^n$-uniformly distributed with constant $2$, and $N\lesssim  (n^{-1}\log(1/\delta))^n$. By \eqref{eqn:heavy incidence assumption}, we have
    \[|\mathcal{I}(\mathcal{P},\mathcal{T}_i)|\geqslant \delta^{1-\kappa}|\mathcal{P}||\mathcal{T}_i|.\]

    Since $E$ is sub-uniformly distributed, it is clear that $\mathcal{P}$ is sub-$\{\Delta_j\}_{j=0}^n$-uniform. By Proposition \ref{prop:DOW incidence}, for every $i\in\{1,\dots,N\}$, there exists $\Delta(i)\in[\delta,\delta^\eta]$ such that
    \begin{equation}\label{eqn:incidence of Delta i}
        |\mathcal{I}(\mathcal{P}_{\Delta(i)},\mathcal{T}_{i,\Delta(i)})|\lesssim_\eta \delta^{\kappa-2\eta}\Delta(i)^{-2},
    \end{equation}
    and
    \begin{equation}\label{eqn:density of Delta i}
        \frac{|\mathcal{I}(\mathcal{P}_{\Delta(i)},\mathcal{T}_{i,\Delta(i)})|}{\Delta(i)|\mathcal{P}_{\Delta(i)}||\mathcal{T}_{i,\Delta(i)}|}\gtrsim_\eta \delta^{\eta-\kappa}.
    \end{equation}

    Note that a standard tube with width $10\Delta(i)$ meets at most $O(\Delta(i)^{-1})$ dyadic $\Delta(i)$-cubes in $\mathcal{D}_{\Delta(i)}$. Also, if a disc $B_\mathbf p$, that has the same center as $\mathbf p\in \mathcal{D}_{\Delta(i)}$ with diameter $10\Delta(i)$, meets a standard tube with width $10\Delta(i)$, then $\mathbf p$ meets an $O(\Delta(i))$-neighborhood of that tube, which also contains at most $O(\Delta(i)^{-1})$ dyadic $\Delta(i)$-cubes. Thus, we have
    \[|\mathcal{I}(\mathcal{P}_{\Delta(i)},\mathcal{T}_{i,\Delta(i)})|\lesssim \Delta(i)^{-1}|\mathcal{T}_{i,\Delta(i)}|,\]
    so
    \[\frac{|\mathcal{I}(\mathcal{P}_{\Delta(i)},\mathcal{T}_{i,\Delta(i)})|}{\Delta(i)|\mathcal{P}_{\Delta(i)}||\mathcal{T}_{i,\Delta(i)})|}\leqslant \Delta(i)^{-2}|\mathcal{P}_{\Delta(i)}|^{-1}\lesssim_v \Delta(i)^{-2+s-v}.\]
    Here we use $|\mathcal{P}_{\Delta(i)}|=|E|_{\Delta(i)}\gtrsim_v \Delta(i)^{-s+v}$.
    Combining with \eqref{eqn:density of Delta i}, we have
    \begin{equation}\label{eqn:bound on Delta i}
        \Delta(i)\lesssim_{\eta,v} \delta^{\frac{\kappa-\eta}{2-s+v}}.
    \end{equation}

    For each family $\mathcal{T}_i$, define
    \[E_{\theta,\mathcal{T}_i}:=E\cap \bigcup_{T\in\mathcal{T}_i:\theta\in\Theta_T} \bbT_T.\]
    Then 
    \[E_{\theta,\mathcal{T}}=\bigcup_{i=1}^N E_{\theta,\mathcal{T}_i}.\]
    For each dyadic tube $\mathbf T\in \mathcal{T}_{i,\Delta(i)}$, define $\Theta_\mathbf T$ as its direction set, similar to $\Theta_T$ for $T\in \mathcal{T}$. Note that we have $\mathcal{H}^1(\Theta_\mathbf T)\lesssim \Delta(i)$.

    We claim that
    \[E_{\theta,\mathcal{T}_i}\subseteq \bigcup\left\{\mathbf p\in \mathcal{P}_{\Delta(i)}:\exists\, \mathbf T\in \mathcal{T}_{i,\Delta(i)}\ \text{such that}\ (\mathbf p,\mathbf T)\in \mathcal{I}(\mathcal{P}_{\Delta(i)}, \mathcal{T}_{i,\Delta(i)}),\ \theta\in \Theta_\mathbf T\right\}.\]
    Indeed, let $x\in E_{\theta,\mathcal{T}_i}$. Then $x\in E\cap \bbT_T$ and $\theta\in\Theta_T$ for some $T\in\mathcal{T}_i$. Let $p\in\mathcal{P}$ be the unique $\delta$-cube containing $x$. Then $x\in B_p\cap\bbT_T$, so $(p,T)\in \mathcal{I}(\mathcal{P},\mathcal{T}_i)$. Let $(\mathbf{p},\mathbf T)$ be the $\Delta(i)$-parents of $(p,T)$. Then $x\in \mathbf p$. By Lemma \ref{lem:containment}, we have $(\mathcal{I}(\mathcal{P},\mathcal{T}_{i}))_{\Delta(i)}\subseteq \mathcal{I}(\mathcal{P}_{\Delta(i)},\mathcal{T}_{i,\Delta(i)})$, so $(\mathbf{p},\mathbf T)\in \mathcal{I}(\mathcal{P}_{\Delta(i)},\mathcal{T}_{i,\Delta(i)})$ with $\theta\in \Theta_T\subseteq\Theta_\mathbf T$.

    Thus, we have
    \[\mathcal{H}^t_\infty(E_{\theta, \mathcal{T}_i})\lesssim \Delta(i)^t\sum_{(\mathbf p,\mathbf T)\in \mathcal{I}(\mathcal{P}_{\Delta(i)},\mathcal{T}_{i,\Delta(i)})} 1_{\Theta_\mathbf{T}}(\theta).\]
    By \eqref{eqn:incidence of Delta i}, it follows that
    \begin{align*}
        \int_{\mathbb S^1} \mathcal{H}^t_\infty(E_{\theta, \mathcal{T}_i})\,d\mathcal{H}^1(\theta)&\lesssim\Delta(i)^{t}\int_{\mathbb S^1}\sum_{(\mathbf p,\mathbf T)\in \mathcal{I}(\mathcal{P}_{\Delta(i)},\mathcal{T}_{i,\Delta(i)})} 1_{\Theta_\mathbf{T}}(\theta)\,d\mathcal{H}^1(\theta)\\
        &=\Delta(i)^t \sum_{(\mathbf p,\mathbf T)\in \mathcal{I}(\mathcal{P}_{\Delta(i)},\mathcal{T}_{i,\Delta(i)})} \mathcal{H}^1(\Theta_{\mathbf{T}})\\
        &\lesssim \Delta(i)^{t+1}|\mathcal{I}(\mathcal{P}_{\Delta(i)},\mathcal{T}_{i,\Delta(i)})|\\
        &\lesssim_\eta \delta^{\kappa-2\eta}\Delta(i)^{t-1}.
    \end{align*}
    Since $t>1$, by \eqref{eqn:bound on Delta i} we have
    \[\Delta(i)^{t-1}\lesssim_{\eta,v} \delta^{\frac{(t-1)(\kappa-\eta)}{2-s+v}},\]
    so
    \begin{align*}
        \int_{\mathbb S^1} \mathcal{H}^t_\infty(E_{\theta, \mathcal{T}})\,d\mathcal{H}^1(\theta)&\leqslant\sum_{i=1}^N\int_{\mathbb S^1} \mathcal{H}^t_\infty(E_{\theta, \mathcal{T}_i})\,d\mathcal{H}^1(\theta)\\
        &\lesssim_{\eta,v} N\cdot\delta^{\kappa-2\eta+\frac{(t-1)(\kappa-\eta)}{2-s+v}}\\
        &=N\cdot \delta^{\kappa\left(1+\frac{t-1}{2-s+v}\right)-\eta\left(2+\frac{t-1}{2-s+v}\right)}.
        %&\lesssim \delta^{\kappa-\eta+\frac{(t-1)(\kappa-\eta)}{2-s+v}}.
    \end{align*}
    %The final step we use $N\lesssim (n^{-1}\log(1/\delta))^{n}\leqslant \log(1/\delta)^{n}\lesssim \delta^{-\eta}$.
    Note that
    \[\kappa\left(1+\frac{t-1}{2-s+v}\right)-\eta\left(2+\frac{t-1}{2-s+v}\right)>\kappa\left(1+\frac{t-1}{2-s}\right)-\frac{2\varepsilon}{3}=\frac{\kappa(t+1-s)}{2-s}-\frac{2\varepsilon}{3}.\]
    Also, $N\lesssim (n^{-1}\log(1/\delta))^{n}\leqslant \log(1/\delta)^{n}\lesssim_{\varepsilon,\eta} \delta^{-\varepsilon/3}$.
    Therefore,
    \[\int_{\mathbb S^1} \mathcal{H}^t_\infty(E_{\theta, \mathcal{T}})\,d\mathcal{H}^1(\theta)\lesssim_{s,t,\kappa,\varepsilon} \delta^{\frac{\kappa(t+1-s)}{2-s}-\varepsilon}.\qedhere\]
\end{proof}

\begin{rmk}
    One may think that using the inequality $|\mathcal{I}(\mathcal{P}_{\Delta(i)},\mathcal{T}_{i,\Delta(i)})|\lesssim \Delta(i)^{-1}|\mathcal{T}_{i,\Delta(i)}|$ to estimate $\Delta(i)$ is too crude. In fact, combining with \eqref{eqn:density of Delta i}, we can use ``high frequency term'' in the proof of Proposition \ref{prop:DOW incidence} to get
    \[|\mathcal{I}(\mathcal{P}_{\Delta(i)},\mathcal{T}_{i,\Delta(i)})|\lesssim_\eta \sqrt{S\Delta(i)^{-1}|\mathcal{P}_{\Delta(i)}||\mathcal{T}_{i,\Delta(i)}}|.\]
    Thus, we have
    \[\Delta(i)^{3-s+v}|\mathcal{T}_{i,\Delta(i)}|\lesssim_{\eta,v}\delta^{2\kappa-3\eta}.\]
    Using the trivial bound $|\mathcal{T}_{i,\Delta(i)}|\geqslant 1$, we have
    \[\Delta(i)\lesssim_{\eta,v}\delta^{\frac{2\kappa-3\eta}{3-s+v}},\]
    which is weaker than $\Delta(i)\lesssim_{\eta,v} \delta^{\frac{\kappa-\eta}{2-s+v}}$. If one can obtain a lower bound on $|\mathcal{T}_{i,\Delta(i)}|$, then this estimate can be improved.
\end{rmk}

\section{Proof of Theorem \ref{thm:main result}}\label{sec:proof of main}
\subsection{Setup}\label{setup}
By scaling and translation, we may assume that $E\subseteq[0,1)^2$ and $\mu(E)\sim 1$. We also assume $1<s<2$, because the case $s=1$ is trivial and for case $s=2$, $s-\alpha_s(s-1)=9/5>3/2$.

Let $0<\varepsilon\ll 1/1000$, $0<\alpha<1$ (we later fix $\alpha=1/(3+\sqrt{8-2s})$). Define
\[\tau:=\alpha(s-1)-5\varepsilon.\]
We will show that
\begin{equation}\label{eqn:Haus content 0}
    \mathcal{H}^{s-\tau}_\infty (\Vis_\theta(E))=0,\quad\text{for $\mathcal{H}^1$-a.e.}\ \theta\in \mathbb S^1
\end{equation}
Taking $\varepsilon\to 0$ gives \eqref{eqn:main bound}.

Let $\delta\in 2^{-\mathbb N}$ be sufficiently small. Let $0<h,k<1$ that will be fixed later. Choose $\Delta\in 2^{-\mathbb N}$ with $\delta^k\leqslant \Delta \leqslant 2\delta^k$. Note that $\delta\ll \Delta$.

Since $E$ is $s$-Ahlfors-regular, we have $|\mathcal{D}_\delta(E)|\sim \delta^{-s}$ and $|\mathcal{D}_\Delta(E)|\sim \delta^{-ks}$. For every $Q\in\mathcal{D}_\Delta(E)$, define the measure $\mu_Q:=\mu|_{3Q}$. Note that $\mu_Q(\mathbb R^2)=\mu(3Q)\sim \delta^{ks}$.

By Proposition \ref{prop:sobolev proj}, we have
\begin{equation}\label{eqn:sobolev estimate global}
    \int_{\mathbb S^1} \|\mu_{\theta^\perp}\|_{H^{(s-1-\varepsilon)/2}(\mathbb R)}^2\,d\mathcal{H}^1(\theta)\lesssim_\varepsilon 1,
\end{equation}
and
\begin{equation}\label{eqn:sobolev estimate local}
    \int_{\mathbb S^1} \|\mu_{Q,\theta^\perp}\|_{H^{(s-1-\varepsilon)/2}(\mathbb R^2)}^2\,d\mathcal{H}^1(\theta)\lesssim_\varepsilon\mu_Q(\mathbb R)\sim \delta^{ks},\quad\forall\,Q\in\mathcal{D}_\Delta(E).
\end{equation}

Now we consider the standard tubes at direction $\theta$.
\begin{notat*}[Standard tubes at direction $\theta$]{}
    Given $\theta\in\mathbb S^1$ and $\gamma\in 2^{-\mathbb N}$, let the collection of standard $\gamma$-tubes at the direction $\theta$ be
    \[\mathscr{T}_{\theta,\gamma}:=\{\bbT=P_{\theta^\perp}^{-1}(2I):I\in \mathcal{D}_\gamma(\mathbb R),\bbT\cap [0,1]^2\ne\varnothing\}.\]
    Note that $E\subseteq [0,1]^2\subseteq \bigcup\limits_{\bbT\in \mathscr{T}_{\theta,\gamma}} \bbT$ and $|\mathscr{T}_{\theta,\gamma}|\sim \gamma^{-1}$. Also, let $\mathscr{T}_{\theta}:=\bigcup\limits_{\gamma}\mathscr{T}_{\theta,\gamma}$.
\end{notat*}

By classifying the standard tubes at direction $\theta$, we will first decompose $\Vis_\theta(E)$ into light part, heavy part, good part, and bad part.

\subsection{Classification of light tubes and heavy tubes}
We first classify light tubes and heavy tubes.

\subsubsection{Light tubes}
Let the collection of light tubes at direction $\theta$ be
\[\mathscr{T}_{\theta,\delta,L}:=\{\bbT\in\mathscr{T}_{\theta,\delta}: |\bbT\cap E|_\delta\leqslant\delta^{-s+1+\tau+\varepsilon}\}.\]
Let the light part of $E$ associated with tubes at direction $\theta$ be 
\[E_{\theta,L}:=E\cap\bigcup_{\bbT\in\mathscr{T}_{\theta,\delta,L}} \bbT.\]

\begin{lem}[Light part estimate]{}\label{lem:light part}
    For every $\theta\in \mathbb S^1$, we have
    \[\mathcal{H}^{s-\tau}_\infty(E_{\theta,L})\lesssim\delta^\varepsilon.\]
\end{lem}

\begin{proof}
    Note that
    \[|E_{\theta,L}|_\delta\leqslant \sum_{T\in\mathscr{T}_{\theta,\delta,L}} |\bbT\cap E|_\delta\lesssim \delta^{-1}\cdot\delta^{-s+1+\tau+\varepsilon}=\delta^{-s+\tau+\varepsilon}.\]
    Then
    \[\mathcal{H}^{s-\tau}_\infty(E_{\theta,L})\lesssim |E_{\theta,L}|_\delta\cdot\delta^{s-\tau}=\delta^\varepsilon.\qedhere\]
\end{proof}

\subsubsection{Heavy tubes and heavy part}\label{heavy set-up}
Since $E$ is $s$-Ahlfors-regular set, we have $|E|_\gamma\sim \gamma^{-s}$, so we expect that for a typical tube $\bbT\in \mathscr{T}_{\theta,\gamma}$, we have
\[|\bbT\cap E|_\gamma\lesssim \gamma^{-s+1}.\]
We classify the heavy tubes whose covering number is very large compared to the above expectation.

Let the collection of heavy tubes be
\[\mathscr{T}_{\theta,H}:=\{\bbT\in\mathscr{T}_\theta:|\bbT\cap E|_{w(\bbT)}\geqslant \delta^{-2\varepsilon}w(\bbT)^{-s+1}\}.\]
Let the heavy part of $E$ associated with tubes at direction $\theta$ be
\[E_{\theta,H,\mathscr{T}}:=E\cap \bigcup_{\bbT\in\mathscr{T}_{\theta,H}} \bbT\]

For every $Q\in\mathcal{D}_\Delta(E)$, let the collection of heavy tubes inside $Q$ be
\[\mathscr{T}_{\theta,\delta,H}(Q):=\{\bbT\in \mathscr{T}_{\theta,\delta}:|\bbT\cap E\cap Q|_\delta\geqslant \delta^{(k-1)(s-1)-h}\}.\]
Let the heavy part of $Q$ associated with tubes at direction $\theta$ be
\[Q_{\theta,H}=E\cap 3Q\cap \bigcup_{\bbT\in\mathscr{T}_{\theta,\delta,H}(Q)} \bbT.\]

Finally, let the heavy part of $E$ at direction $\theta$ be
\[E_{\theta,H}:=E_{\theta,H,\mathscr{T}}\cup \bigcup_{Q\in\mathcal{D}_\Delta(E)} Q_{\theta,H}.\]

\subsection{Decomposition of good part and bad part}
For every $Q\in \mathcal{D}_\Delta(E)$, let the collection of rich tubes associated with $Q$ be
\[\mathscr{T}_{\theta,\delta}(Q):=\left\{\bbT\in\mathscr{T}_{\theta,\delta}:|\bbT\cap E\cap Q|_\delta\geqslant \delta^{(k-1)(s-1)+\tau+4\varepsilon}\right\}.\]

\begin{lem}[Rich tube]{}\label{lem:rich tube}
    Let $\bbT\in\mathscr{T}_{\theta,\delta}\setminus\mathscr{T}_{\theta,\delta,L}$. Suppose $\bbT$ is not contained in any tubes in $\mathscr{T}_{\theta,H}$. Then there exists $Q\in\mathcal{D}_\Delta(E)$ such that $\bbT\in\mathscr{T}_{\theta,\delta}(Q)$.
\end{lem}

\begin{proof}
    Let $\bbT'\in \mathscr{T}_{\theta,\Delta}$ with $\bbT\subseteq \bbT'$. Since $\bbT'\notin\mathscr{T}_{\theta,H}$, we have
    \[|\bbT\cap E|_\Delta\leqslant |\bbT'\cap E|_\Delta\lesssim\delta^{-2\varepsilon}\Delta^{-s+1}\sim \delta^{-k(s-1)-2\varepsilon}.\]
    Since $\bbT\notin\mathscr{T}_{\theta,\delta,L}$, we have $|\bbT\cap E|_\delta>\delta^{-s+1+\tau+\varepsilon}$. By the pigeonhole principle, there exists $Q\in\mathcal{D}_\Delta(E)$ such that
    \[|\bbT\cap E\cap Q|_\delta\gtrsim\frac{|\bbT\cap E|_\delta}{|\bbT\cap E|_\Delta}\gtrsim \delta^{(k-1)(s-1)+\tau+3\varepsilon}.\qedhere\]
\end{proof}

Let the collection of all affine lines at direction $\theta$ be $\mathscr{L}_\theta$.
For every $Q\in\mathcal{D}_\Delta(E)$, let $\mathscr{L}_{\theta,Q,B}$ be the collection of bad lines $\ell\in\mathscr{L}_\theta$ such that there exists a $\bbT\in\mathscr{T}_{\theta,\delta}(Q)$ with the following properties:
\begin{enumerate}[(1)]
    \item $\ell\subseteq \bbT$,
    \item $\ell\cap \overline{3Q}\cap E=\varnothing$,
\end{enumerate}
Set
\[L_{\theta,Q,B}:=\bigcup_{\ell\in \mathscr{L}_{\theta,Q,B}}\ell\quad\text{and}\quad L_{\theta,B}:=\bigcup_{Q\in\mathcal{D}_\Delta(E)} L_{\theta,Q,B},\]
and let the bad part of $E$ at direction $\theta$ be
\[E_{\theta,B}:=E\cap L_{\theta,B}\setminus E_{\theta,H}.\]

Let the good part of $E$ at direction $\theta$ be 
\[E_{\theta,G}:=E\setminus (E_{\theta,H}\cup E_{\theta,L}\cup E_{\theta,B}).\]

Thus, we have
\[\Vis_\theta(E)\subseteq E_{\theta,L}\cup E_{\theta,H}\cup (\Vis_\theta(E)\cap E_{\theta,G})\cup E_{\theta,B}.\]
The Hausdorff content of each part will be estimated in the following sections.

\subsection{Good part}
In this section, we show that the Hausdorff dimension of good part is small.

Let $\bbT\in\mathscr{T}_{\theta,\delta}$ with $\bbT\cap E_{\theta,G}\ne\varnothing$. Then $\bbT\notin \mathscr{T}_{\theta,\delta,L}\cup\mathscr{T}_{\theta,H}$ and $\bbT$ is not contained in any tubes in $\mathscr{T}_{\theta,H}$. By Lemma \ref{lem:rich tube}, there exists a cube $Q\in \mathcal{D}_\Delta(E)$ such that $\bbT\in \mathscr{T}_{\theta,\delta}(Q)$. Abbreviate $\mathcal{Q}_{\Delta}(\bbT):=\mathcal{D}_\Delta(E\cap\bbT)$. Let $\bbT'\in\mathscr{T}_{\theta,\Delta}$ be a tube containing $\bbT$. Then $\bbT'\notin \mathscr{T}_{\theta,H}$, so
\[|\mathcal{Q}_{\Delta}(\bbT)|\lesssim |\bbT'\cap E|_\Delta\lesssim \delta^{-k(s-1)-2\varepsilon}.\]

Let $Q_{\bbT}\in\mathcal{Q}_{\Delta}(\bbT)$ be the $\theta$-highest cube such that $\bbT\in \mathscr{T}_{\theta,\delta}(Q)$ and $\inf\{x\cdot \theta:x\in E\cap 3Q_{\bbT}\}$ is maximized.

\begin{lem}[Visible part is in high cube]\label{lem:Vis not in low cube}
    Let $Q\in\mathcal{Q}_\Delta(\bbT)$. If $Q$ is lower than $Q_{\bbT}$, that is, $\sup\{x\cdot\theta:x\in E\cap 3Q\}<\inf\{x\cdot \theta:x\in E\cap 3Q_{\bbT}\}$, then
    \[\Vis_\theta(E)\cap E_{\theta,G}\cap \bbT\cap Q=\varnothing.\]
\end{lem}

\begin{proof}
    Assume there exists $x\in \Vis_\theta(E)\cap  E_{\theta,G}\cap \bbT\cap Q$. Consider the affine line $\ell\subseteq \bbT$ with $x\in \ell$. Note that $x\notin E_{\theta,B}$, so $x\notin L_{\theta,Q_{\bbT},B}$ and $\ell\notin \mathscr{L}_{\theta,Q_{\bbT},B}$. Since $\bbT\in\mathscr{T}_{\theta,\delta}(Q)$, we have $\ell\cap E\cap \overline{3Q_{\bbT}}\ne\varnothing$. Let $y\in \ell\cap E\cap \overline{3Q_{\bbT}}$. Then $P_{\theta^\perp}(x)=P_{\theta^\perp}(y)$. Since $y\in E\cap 3Q_{\bbT}$, we have $x\cdot\theta<y\cdot \theta$, which leads to a contradiction because $x\in\Vis_\theta(E)$.
\end{proof}

\begin{prop}[$\delta$-covering number of good part]{}\label{prop:delta covering of good part}
    For every $\theta\in \mathbb S^1$, we have
    \[|\Vis_\theta(E)\cap E_{\theta,G} \cap \bbT|_\delta\lesssim \delta^{(k-1)(s-1)-h}+\delta^{(\alpha-1)(s-1)-3\varepsilon}.\]
\end{prop}

\begin{proof}
    By Lemma \ref{lem:Vis not in low cube}, $\Vis_\theta(E)\cap  E_{\theta,G}\cap \bbT$ is in $Q\in\mathcal{Q}_\Delta(\bbT)$ with
    \[\sup\{x\cdot\theta:x\in 3Q\cap E\}\geqslant \inf\{x\cdot\theta:x\in E\cap 3Q_{\bbT}\}.\]
    We split into two cases:
    \begin{enumerate}[(1)]
        \item $\inf\{x\cdot\theta:x\in E\cap3Q\}\leqslant\inf\{x\cdot \theta:x\in E\cap 3Q_{\bbT}\}$,
        \item $\inf\{x\cdot\theta:x\in E\cap3Q\}>\inf \{x\cdot\theta:x\in E\cap 3Q_{\bbT}\}$.
    \end{enumerate}

   In case (1), there exists $\lesssim 1$ cubes satisfying this case, because the number of cubes $Q\in\mathcal{Q}_{\Delta}(\bbT)$ satisfying both inequalities is finite and uniformly bounded. For each such $Q$, we have $E_{\theta,G}\cap Q_{\theta,H}=\varnothing$. If $\bbT\in\mathscr{T}_{\theta,\delta,H}(Q)$, then $E\cap Q\cap \bbT\subseteq Q_{\theta,H}\subseteq E_{\theta,H}$, so $E_{\theta,G}\cap Q\cap \bbT=\varnothing$. If $\bbT\notin\mathscr{T}_{\theta,\delta,H}(Q)$, then
   \[|\Vis_\theta(E)\cap \bbT\cap E_{\theta,G}\cap Q|_\delta\leqslant |E\cap Q\cap \bbT|_\delta\leqslant \delta^{(k-1)(s-1)-h}.\]

   In case (2),
   we have 
   \[|\Vis_\theta(E)\cap  Q \cap \bbT|_\delta\leqslant |E\cap Q\cap \bbT|_\delta\leqslant \delta^{(k-1)(s-1)+\tau+4\varepsilon}.\]
   
   Thus,
   \begin{align*}
       |\Vis_\theta(E)\cap E_{\theta,G} \cap \bbT|_\delta&\lesssim \delta^{(k-1)(s-1)-h}+|\mathcal{Q}_{\Delta}(\bbT)|\cdot \delta^{(k-1)(s-1)+\tau+4\varepsilon}\\
       &\lesssim\delta^{(k-1)(s-1)-h}+\delta^{(\alpha-1)(s-1)-3\varepsilon}.
   \end{align*}
\end{proof}

\begin{cor}[Good part estimate]{}\label{cor:good part}
    For every $\theta\in\mathbb S^1$, we have
    \[\mathcal{H}^{s-\tau}_\infty(\Vis_\theta(E)\cap E_{\theta,G})\lesssim\delta^{(k-\alpha)(s-1)-h+5\varepsilon}+\delta^{2\varepsilon}.\]
\end{cor}

\begin{proof}
    By Proposition \ref{prop:delta covering of good part}, we have
    \begin{align*}
        \mathcal{H}^{s-\tau}_\infty(\Vis_\theta(E)\cap E_{\theta,G})&\lesssim\sum_{\bbT\in\mathscr{T}_{\theta,\delta}}  |\Vis_\theta(E)\cap E_{\theta,G} \cap \bbT|_\delta\cdot\delta^{s-\tau}\\
        &\lesssim \delta^{-1}\cdot (\delta^{(k-1)(s-1)-h}+\delta^{(\alpha-1)(s-1)-3\varepsilon})\cdot\delta^{s-\tau}\\
        &=\delta^{k(s-1)-\tau-h}+\delta^{2\varepsilon}.\\
        &=\delta^{(k-\alpha)(s-1)-h+5\varepsilon}+\delta^{2\varepsilon}.
    \end{align*}
\end{proof}

\subsection{Bad part} 
In this section, we show that the Hausdorff dimension of bad part is small.
\begin{prop}[Bad lines estimate]{}\label{prop:bad lines estimate}
    For every $\theta\in\mathbb S^1$ and $Q\in\mathcal{D}_\Delta(E)$, we have
    \[\mathcal{H}^{1-\tau}_\infty(P_{\theta^\perp}(L_{\theta,Q,B}))\lesssim \delta^{(1-2k-3\alpha)(s-1)+5\varepsilon}\|\mu_{Q,\theta^\perp}\|_{H^{(s-1-\varepsilon)/2}(\mathbb R)}^2.\]
\end{prop}

\begin{proof}
    Suppose $\mathcal{H}^{1-\tau}_\infty(P_{\theta^\perp}(L_{\theta,Q,B}))>0$. Note that $P_{\theta^\perp}(L_{\theta,Q,B})$ is a Borel set. By the Frostman's lemma, there exists a Borel probability measure $\nu$ supported on $P_{\theta^\perp}(L_{\theta,Q,B})$ such that
    \[\nu(B(x,r))\lesssim \mathcal{H}^{1-\tau}_\infty(P_{\theta^\perp}(L_{\theta,Q,B}))^{-1}r^{1-\tau}.\]

    By the definition of bad lines, we have
    \[P_{\theta^\perp}(L_{\theta,Q,B})\cap P_{\theta^\perp}(E\cap\overline{3Q})=\varnothing,\]
    so $\supp\nu\cap \supp \mu_{Q,\theta^\perp}=\varnothing$. Let $\varphi\in C_c^\infty(\mathbb R^n)$ be a standard smooth bump function with $\|\varphi\|_{L^1(\mathbb R)}=1$. Write $\varphi_r(t):=r^{-1}\varphi(t/r)$. Since $\supp\nu$ and $\supp \mu_{Q,\theta^\perp}$ are compact, there exists a small $\eta>0$ so that
    \[\int_{\mathbb R} \varphi_{\eta}* \mu_{Q,\theta^\perp}\,d\nu=0.\]
    Decompose
    \[\varphi_{\eta}*\mu_{Q,\theta^\perp}=(\varphi_{\eta}*\varphi_{C\delta}*\mu_{Q,\theta^\perp})+(\varphi_{\eta}* \mu_{Q,\theta^\perp}-\varphi_{\eta}*\varphi_{C\delta}*\mu_{Q,\theta^\perp}),\]
    where $C\geqslant 1$ is a constant to be fixed later. Let
    \[I_1:=\int_{\mathbb R} \varphi_{\eta}*\varphi_{C\delta}* \mu_{Q,\theta^\perp}\,d\nu\ \text{and}\ I_2:=\int_{\mathbb R} \varphi_{\eta}* \mu_{Q,\theta^\perp}-\varphi_{\eta}*\varphi_{C\delta}*\mu_{Q,\theta^\perp}\,d\nu.\]
    Then $I_1+I_2=0$.

     We first estimate $I_1$. By the definition of bad lines, for every $x\in P_{\theta^\perp}(L_{\theta,Q,B})$, there exists a tube $\bbT\in \mathscr{T}_{\theta,\delta}(Q)$ such that $P_{\theta^\perp}^{-1}(x)\subseteq\bbT$. Thus, we can choose $C\sim 1$ sufficiently large so that, for every $x\in P_{\theta^\perp}(L_{\theta,Q,B})$,
    \[|\ell_x(C\delta)\cap E\cap Q|_\delta\gtrsim \delta^{(k-1)(s-1)+\tau+4\varepsilon},\]
    where $\ell_x:=P_{\theta^\perp}^{-1}(x)$ and $\ell_x(C\delta)$ is the $C\delta$-neighborhood of $\ell_x$. Since $\mu_Q(B(y,\delta))\sim \delta^s$ for $y\in E\cap Q$, we have
    \[\mu_Q(\ell_x(C\delta))\gtrsim |\ell_x(C\delta)\cap E\cap Q|_\delta\cdot\delta^s,\]
    so for every $x\in\supp\nu$, we have
    \[\varphi_\eta*\varphi_{C\delta}*\mu_Q(x)\gtrsim \delta^{-1}\mu_{Q,\theta^\perp}(\ell_x(C\delta))\gtrsim \delta^{k(s-1)+\tau+4\varepsilon}.\]
    Since $\nu(\mathbb R)=1$, we have
    \[I_1\gtrsim \delta^{k(s-1)+\tau+4\varepsilon}.\]

    Now we estimate $I_2$. By the Plancherel's identity, we have
    \[I_2=\int_{\mathbb R} \widehat{\varphi}(\eta\xi)(1-\widehat{\varphi}(C\delta\xi))\widehat{\mu_{Q,\theta^\perp}}(\xi)\overline{\widehat{\nu}}(\xi)\,d\xi.\]
    Since $\widehat{\varphi}$ is a Schwartz function, we have
    \[|1-\widehat{\varphi}(C\delta\xi)|=|\widehat{\varphi}(0)-\widehat{\varphi}(C\delta\xi)|\lesssim\min \{1,\delta|\xi|\}\leqslant \delta^a|\xi|^a,\quad\forall\,0\leqslant a\leqslant 1.\]
    By the Cauchy-Schwarz inequality, we have
    \[|I_2|\lesssim \delta^a\left(\int_{\mathbb R} |\widehat{\mu_{Q,\theta^\perp}}(\xi)|^2|\xi|^{s-1-\varepsilon}\,d\xi\right)^{1/2}\left(\int_{\mathbb R} |\widehat{\nu}(\xi)|^2|\xi|^{2a-s+1+\varepsilon}\,d\xi\right)^{1/2}.\]
    Choose $a=\frac{s-1-\tau}{2}-\varepsilon$. Then we have
    \[\int_{\mathbb R} |\widehat{\nu}(\xi)|^2|\xi|^{-\tau-\varepsilon}\,d\xi\sim I_{1-\tau-\varepsilon}(\nu)\lesssim_\varepsilon\mathcal{H}^{1-\tau}_\infty(P_{\theta^\perp}(L_{\theta,Q,B}))^{-1}, \]
    so
    \[I_2\lesssim_\varepsilon\delta^a\|\mu_{Q,\theta^\perp}\|_{H^{(s-1-\varepsilon)/2}(\mathbb R)}\mathcal{H}^{1-\tau}_\infty(P_{\theta^\perp}(L_{\theta,Q,B}))^{-1/2}.\]

    Since $I_1\leqslant |I_2|$, we have
    \[\delta^{k(s-1)+\tau+4\varepsilon}\lesssim\delta^a\|\mu_{Q,\theta^\perp}\|_{H^{(s-1-\varepsilon)/2}(\mathbb R)}\mathcal{H}^{1-\tau}_\infty(P_{\theta^\perp}(L_{\theta,Q,B}))^{-1/2}.\]
    Simplifying this inequality gives
    \begin{align*}
        \mathcal{H}^{1-\tau}_\infty(P_{\theta^\perp}(L_{\theta,Q,B}))&\lesssim_\varepsilon \delta^{(1-2k)(s-1)-3\tau-10\varepsilon}\|\mu_{Q,\theta^\perp}\|_{H^{(s-1-\varepsilon)/2}(\mathbb R)}^2\\
        &=\delta^{(1-2k-3\alpha)(s-1)+5\varepsilon}\|\mu_{Q,\theta^\perp}\|_{H^{(s-1-\varepsilon)/2}(\mathbb R)}^2.
    \end{align*}
\end{proof}

\begin{cor}[Bad part estimate]{}\label{cor:bad part}
    We have
    \[\int_{\mathbb S^1} \mathcal{H}^{s-\tau}_\infty(E_{\theta,B})\,d\mathcal{H}^1(\theta)\lesssim_\varepsilon  \delta^{(1-2k-3\alpha)(s-1)+3\varepsilon}.\]
\end{cor}

\begin{proof}
    Fix $\theta\in\mathbb S^1$. Let $\mathscr{I}\subseteq\bigcup\limits_{\gamma}\mathcal{D}_\gamma(\mathbb R)$ be a family of dyadic intervals in $\mathbb R$ such that $P_{\theta^\perp}(L_{\theta,B})\subseteq \bigcup\limits_{I\subseteq \mathscr{I}} I$ and
    \begin{equation}\label{eqn:cover of proj of bad part}
        \sum_{I\in\mathscr{I}}|I|^{1-\tau}\lesssim\sum_{Q\in\mathcal{D}_\Delta(E)}\mathcal{H}^{1-\tau}_\infty(P_{\theta^\perp}(L_{\theta,B})).
    \end{equation}
    Since $P_{\theta^\perp}(E_{\theta,B})\subseteq P_{\theta^\perp}(L_{\theta,B})$, we may assume $I\cap P_{\theta^\perp}(E_{\theta,B})\ne\varnothing$ for every $I\in \mathscr{I}$.

    For each $I\in\mathscr{I}$, let $\bbT_I\in\mathcal{T}_\theta$ with $\bbT_I=P_{\theta^\perp}^{-1}(2I)$. Since $I\cap P_{\theta^\perp}(E_{\theta,B})\ne\varnothing$, we have $\bbT_I\cap E_{\theta,B}\ne\varnothing$. Note that $E_{\theta,B}\subseteq E\setminus E_{\theta,H}$, so $\bbT_I\notin \mathscr{T}_{\theta,H}$. Thus,
    \begin{equation}\label{eqn:cover not in heavy part}
        |\bbT_I\cap E|_{|I|}\leqslant \delta^{-2\varepsilon}|I|^{-s+1}.
    \end{equation}
    Since $E_{\theta,B}\subseteq E\cap\bigcup\limits_{I\in\mathscr{I}}\bbT_I$, by \eqref{eqn:cover of proj of bad part} and \eqref{eqn:cover not in heavy part}, we have
    \begin{align*}
        \mathcal{H}^{s-\tau}_\infty(E_{\theta,B})&\lesssim\sum_{I\in\mathscr{I}}|\bbT_I\cap E|_{|I|}\cdot |I|^{s-\tau}\\
        &\leqslant \sum_{I\in\mathscr{I}}\delta^{-2\varepsilon}|I|^{1-\tau}\\
        &\lesssim \delta^{-2\varepsilon}\sum_{Q\in\mathcal{D}_\Delta(E)}\mathcal{H}^{1-\tau}_\infty(P_{\theta^\perp}(L_{\theta,Q,B}))\\
        &\lesssim_\varepsilon \delta^{(1-2k-3\alpha)(s-1)+3\varepsilon}\sum_{Q\in\mathcal{D}_\Delta(E)}\|\mu_{Q,\theta^\perp}\|_{H^{(s-1-\varepsilon)/2}(\mathbb R)}^2.
    \end{align*}

    Therefore, by \eqref{eqn:sobolev estimate local}, we have
    \begin{align*}
        \int_{\mathbb S^1} \mathcal{H}^{s-\tau}_\infty(E_{\theta,B})\,d\mathcal{H}^1(\theta)&\lesssim \delta^{(1-2k-3\alpha)(s-1)+3\varepsilon}\sum_{Q\in\mathcal{D}_\Delta(E)}\int_{\mathbb S^1}\|\mu_{Q,\theta^\perp}\|_{H^{(s-1-\varepsilon)/2}(\mathbb R)}^2 \,d\mathcal{H}^1(\theta)\\
        &\lesssim\delta^{(1-2k-3\alpha)(s-1)+3\varepsilon}\sum_{Q\in\mathcal{D}_\Delta(E)} \mu_Q(\mathbb R^2)\\
        &\sim \delta^{(1-2k-3\alpha)(s-1)+3\varepsilon}.
    \end{align*}
\end{proof}

\subsection{Heavy part}
Recall the notation in \ref{heavy set-up}.
In this section, we estimate $E_{\theta,H,\mathscr{T}}$ and $Q_{\theta,H}$ separately.

\subsubsection{Global heavy part}
We apply the maximal function estimate in Proposition \ref{prop:maximal function estimate} to control the size of $E_{\theta,H,\mathscr{T}}$.

\begin{lem}[Maximal level set]\label{lem:global heavy containment}
    For every $\theta\in\mathbb S^1$, we have
    \[E_{\theta,H,\mathscr{T}}\subseteq \{x\in E:\mathcal{M}\mu_{\theta^\perp}(P_{\theta^\perp}(x))\geqslant \delta^{-\varepsilon}\}.\]
\end{lem}

\begin{proof}
    Let $x\in E_{\theta,H}$. Then there exists $\bbT\in \mathscr{T}_{\theta}$ with $x\in \bbT$. Then there exists $C\sim 1$ such that
    \begin{align*}
        \mathcal{M}\mu_{\theta^\perp}(P_{\theta^\perp}(x))&\geqslant (Cw(\bbT))^{-1}\mu_{\theta^\perp}(B(P_{\theta^\perp}(x),Cw(\bbT)))\\
        &\gtrsim w(\bbT)^{-1}\mu(C\bbT)\\
        &\gtrsim w(\bbT)^{-1}|\bbT\cap E|_{w(\bbT)}\cdot w(\bbT)^s\\
        &\geqslant\delta^{-2\varepsilon}.
    \end{align*}
    Here we use $s$-Ahlfors-regularity of $\mu$.
\end{proof}

\begin{cor}[Global heavy part estimate]\label{cor:global heavy part}
    We have
    \[\int_{\mathbb S^1} \mathcal{H}^{s-\tau}_\infty(E_{\theta,H,\mathscr{T}})\,d\mathcal{H}^1(\theta)\lesssim_\varepsilon\delta^\varepsilon.\]
\end{cor}

\begin{proof}
    Fix $\theta\in \mathbb S^1$. Set $H_\theta(E):=\{x\in E:\mathcal{M}\mu_{\theta^\perp}(P_{\theta^\perp}(x))\geqslant \delta^{-\varepsilon}\}$.
    By Lemma \ref{lem:global heavy containment} and Proposition \ref{prop:maximal function estimate}, we have
    \[\mathcal{H}^{1+2\varepsilon}_\infty(E_{\theta,H,\mathscr{T}})\leqslant \mathcal{H}^{1+2\varepsilon}_\infty(H_\theta(E))\lesssim_\varepsilon\delta^\varepsilon\|\mu_{\theta^\perp}\|_{H^{(s-1-\varepsilon)/2}(\mathbb R)}^2.\]

    Thus, by \eqref{eqn:sobolev estimate global}, we have
    \[\int_{\mathbb S^1} \mathcal{H}^{1+2\varepsilon}_\infty(E_{\theta,H,\mathscr{T}})\,d\mathcal{H}^1(\theta)\lesssim \delta^\varepsilon\int_{\mathbb S^1}\|\mu_{\theta^\perp}\|_{H^{(s-1-\varepsilon)/2}(\mathbb R)}^2\,d\mathcal{H}^1(\theta)\lesssim_\varepsilon \delta^\varepsilon\mu(E)\sim \delta^\varepsilon.\]
    Recall the following inequality
    \[\mathcal{H}^t_\infty(A)\leqslant\diam(A)^{t-r}\mathcal{H}^{r}_\infty(A),\quad\forall\,0<r<t.\]
    Since $1+2\varepsilon<s-\tau$, we have
    \[\int_{\mathbb S^1} \mathcal{H}^{s-\tau}_\infty(E_{\theta,H,\mathscr{T}})\,d\mathcal{H}^1(\theta)\lesssim_\varepsilon\delta^\varepsilon.\qedhere\]
\end{proof}

\subsubsection{Local heavy part}
%In order to apply Proposition \ref{prop:local heavy from incidence} on $Q_{\theta,H}$, we first check $E\cap \overline{3Q}$ is sub-uniformly distributed under affine transformations.
For the local heavy part estimate, we have two scales, $\delta$ and $\Delta$. 
 In order to apply Proposition \ref{prop:local heavy from incidence}, we scale the set $E\cap \overline{3Q}$ back to a subset of the unit cube, and check it is sub-uniformly distributed. 
% For convenience, 

\begin{lem}[Local sub-uniform distribution]{}\label{lem:local scale subuniform}
    Let $Q\in\mathcal{D}_\Delta(E)$, and $\Phi$ be an affine transformation with scaling factor $C\Delta^{-1}$ such that $\Phi(E\cap \overline{3Q})\subseteq [1/8,1/4]\times[1/2,3/4]$. Then for $0<r\leqslant R\leqslant 1$, $|\Phi(E\cap3Q)|_R\sim R^{-s}$ and $|\Phi(E\cap\overline{3Q})\cap p|_r\lesssim (R/r)^s$ for every cube $p$ with side length $R$. 
    Consequently,
    \[|\Phi(E\cap\overline{3Q})|_R\cdot|\Phi(E\cap\overline{3Q})\cap p|_r\lesssim |\Phi(E\cap\overline{3Q})|_r,\]
    so $\Phi(E\cap3Q)$ is sub-uniformly distributed. Moreover, $\dim_\mathrm{H}(\Phi(E\cap\overline{3Q}))=\dim_\mathrm{H}(E\cap3Q)=s$.
\end{lem}

For reasons that will become apparent later on, we scale back to $[1/8,1/4]\times[1/2,3/4]$. 

\begin{proof}
    Abbreviate $E_Q:=\Phi(E\cap \overline{3Q})$. Note that $\mu(\overline{3Q})\sim \Delta^s$ and $|\Phi(x)-\Phi(y)|=C\Delta^{-1}|x-y|$.

    For every $p\in\mathcal{D}_R(E_Q)$, we choose $x_p\in E_Q\cap p$ and let $y_p:=\Phi^{-1}(x_p)\in E$. Then we have $\Phi^{-1}(p)\subseteq B(y_p,2C^{-1}\Delta R)$. Since $\{\Phi^{-1}(p):p\in \mathcal{D}_R(E_Q)\}$ covers $E\cap \overline{3Q}$, we have
    \[\Delta^s\sim \mu(\overline{3Q})\leqslant\sum_{p\in \mathcal{D}_R(E_Q)} \mu(\Phi^{-1}(p))\lesssim \Delta^sR^s|E_Q|_R\implies |E_Q|_R\gtrsim R^{-s}.\]

    Let $0<r\leqslant R\leqslant 1$, and $p$ be a cube with side length $R$ and $E_Q\cap p\ne\varnothing$. Choose a maximal $r$-separated set $\{x_1,\dots,x_N\}\subseteq E_Q\cap p$. Then $|E_Q\cap p|_r\lesssim N$. Let $y_j:=\Phi^{-1}(x_j)\in E$ and $t:=(4C)^{-1}\Delta r$. Then the balls $\{B(y_j,t)\}$ are pairwise disjoint. Moreover, since $x_j\in p$ and $r<R$, we have $\bigcup\limits_{j=1}^N B(y_j,t)\subseteq B(y_1,3C^{-1}\Delta R)$. Thus,
    \[Nt^s\sim \sum_{j=1}^N \mu(B(y_j,t))\leqslant \mu(B(y_1,3C^{-1}\Delta R))\sim \Delta^s R^s \implies N\lesssim (R/r)^s.\]
    Therefore, $|E_Q\cap p|_r\lesssim (R/r)^s$ and $|E_Q|_R\lesssim R^{-s}$.
\end{proof}

Now we are ready to apply Proposition \ref{prop:local heavy from incidence}. The local heavy part assumption transfers nicely to the assumptions of Proposition \ref{prop:local heavy from incidence}.

\begin{prop}[Local heavy part estimate]\label{prop:local heavy part}
    For every sufficiently small $\varepsilon_0>0$, we have
    \[\sum_{Q\in\mathcal{D}_\Delta(E)}\int_{\mathbb S^1} \mathcal{H}^{s-\tau}_\infty(Q_{\theta,H})\, d\mathcal{H}^1(\theta)\lesssim_{\varepsilon,\varepsilon_0} \delta^{\frac{h(1-\alpha s+\alpha)}{2-s}-k\alpha(s-1)+5\left(\frac{ h}{2-s}+k\right)\varepsilon-\varepsilon_0}.\]
\end{prop}

\begin{proof}
    Fix $Q\in\mathcal{D}_\Delta(E)$. In order to apply Proposition \ref{prop:local heavy from incidence}, we first enlarge (translate and scale) $E\cap \overline{3Q}$ to $[1/8,1/4]\times[1/2,3/4]$ by an affine transformation $\Phi_Q$ with scaling factor $C\Delta^{-1}$. Since $E\subseteq[0,1)^2$, we have $C\sim 1$ is uniform for every $Q\in\mathcal{D}_\Delta(E)$. Let $E_Q:=\Phi_Q(E\cap \overline{3Q})$. By Lemma \ref{lem:local scale subuniform}, $E_Q$ is sub-uniformly distributed and $\dim_\mathrm{H} E_Q=s$.

    We only work on the tubes $\bbT\in \bigcup\limits_{\theta\in \mathbb S^1} \mathscr{T}_{\theta,\delta,H}(Q)$ whose directions have slope $a\in [0,1)$. This restriction does not affect our application. The rest of the tubes can be estimated by decomposing angles and applying a rotation. The sub-uniform distributed property and rich tube assumption are preserved under rotations.

    Let $r:=\frac{C\delta}{\Delta}\sim \delta^{1-k}$. Since $\delta\ll\Delta$, we may assume $r$ is sufficiently small. 
  To apply Proposition \ref{prop:local heavy from incidence}, we first convert the standard tubes to dyadic tubes by detecting the center lines.
    For every $\bbT$, let $\bbT_\Phi:=\Phi_Q(\bbT)$. Then $\bbT_\Phi$ is a standard tube with $w(\bbT_\Phi)=Cr$. Suppose the center line of $\bbT_\Phi$ is $y=ax+b$. Since $E_Q\cap \bbT_\Phi\ne\varnothing$, for sufficiently small $r$, we have $0<b<1$. Choose unique $p_{\bbT}=[a_{\bbT},a_{\bbT}+r)\times [b_{\bbT},b_{\bbT}+r)\in\mathcal{D}_r([0,1)^2)$ such that $(a,b)\in p_{\bbT}$. Let the dyadic $r$-tube $T_{\bbT}:=\mathbf D(p_{\bbT})$, and the standard tube associated with it be $\mathbf{T}_{\bbT}$.
    Recall $\mathbf{T}_{\bbT}$ has the center line $y=ax+b$ with width $10r$.
    We claim that
    \[\bbT_\Phi\cap B(0,1)\subseteq \mathbf{T}_{\bbT}.\]
    Indeed, let $(u,v)\in \bbT_\Phi\cap B(0,1)$. Then $|v-au-b|\leqslant (1+a^2)^{1/2}r/2\leqslant \sqrt{2}r/2$. Then
    \[|v-a_{\bbT} u-b_{\bbT}|\leqslant |v-au-b|+|a-a_{\bbT}||u|+|b-b_{\bbT}|\leqslant \sqrt{2}r/2+r+r<3r.\]
    Since $w(\mathbf{T}_{\bbT})=10r$, we have $(u,v)\in \mathbf{T}_{\bbT}$.

    Since $|E\cap Q\cap \bbT|_\delta\geqslant \delta^{(k-1)(s-1)-h}$, we have
    \[|\Phi_Q(E\cap Q)\cap \bbT_\Phi|_r\gtrsim \delta^{(k-1)(s-1)-h}\gtrsim r^{1-s}\delta^{-h}.\]
    Let $x\in p\cap\Phi_Q(E\cap Q)\cap \bbT_\Phi$ for some $p\in \mathcal{D}_r(E_Q)$. Since $\bbT_\Phi\cap B(0,1)\subseteq \mathbf{T}_{\bbT}$, we have $x\in B_p\cap \mathbf{T}_{\bbT}$. By Lemma \ref{lem:local scale subuniform}, we have $|\mathcal{D}_r(E_Q)|=|E_Q|_r\sim r^{-s}$. It follows that
    \begin{align*}
        |\mathcal{I}(\mathcal{D}_r(E_Q),\{T_{\bbT}\})|&=|\{p\in \mathcal{D}_r(E_Q):B_p\cap\mathbf{T}_{\bbT}\neq\varnothing\}|\\
        &\gtrsim r^{1-s}\delta^{-h}\sim r\delta^{-h}|\mathcal{D}_r(E_Q)|\\
        &=r^{1-\frac{h}{1-k}}|\mathcal{D}_r(E_Q)|.
    \end{align*}

    Let $\kappa:=\frac{h-\varepsilon_1}{1-k}$ for sufficiently small $\varepsilon_1>0$ (our final choices of $h$ and $k$ guarantee $\kappa<1$), so that
    \[\mathcal{I}(\mathcal{D}_r(E_Q),\{T_{\bbT}\})\geqslant r^{1-\kappa}|\mathcal{D}_r(E_Q)|.\]
    By Proposition \ref{prop:local heavy from incidence}, for every $\varepsilon_2>0$, we have
    \begin{align*}
        \int_{\mathbb S^1} \mathcal{H}^{s-\tau}_\infty(\Phi_Q(Q_{\theta,H}))\,d\mathcal{H}^1(\theta)&\lesssim r^{\frac{\kappa(1-\tau)}{2-s}-\varepsilon_2}\\
        &=r^{\frac{(h-\varepsilon_1)(1-\tau)}{(1-k)(2-s)}-\varepsilon_2}\\
        &\sim \delta^{\frac{(h-\varepsilon_1)(1-\tau)}{2-s}-(1-k)\varepsilon_2}\\
        &\sim \delta^{\frac{h(1-\tau)}{2-s}-\varepsilon_0},
    \end{align*}
    where
    \[\varepsilon_0:=\frac{\varepsilon_1(1-\tau)}{2-s}+(1-k)\varepsilon_2\leqslant \frac{\varepsilon_1}{2-s}+\varepsilon_2.\]
    We can take $\varepsilon_0>0$ sufficiently small by taking both $\varepsilon_1$ and $\varepsilon_2$ sufficiently small.
    
    Since $\mathcal{H}^{s-\tau}_\infty(\Phi_Q(Q_{\theta,H}))\sim \Delta^{\tau-s}\mathcal{H}^{s-\tau}_\infty(Q_{\theta,H})$, we have
    \[\int_{\mathbb S^1} \mathcal{H}^{s-\tau}_\infty(Q_{\theta,H})\,d\mathcal{H}^1(\theta)\lesssim \Delta^{s-\tau}\delta^{\frac{h(1-\tau)}{2-s}-\varepsilon_0}.\]
    Summing over all $Q\in\mathcal{D}_\Delta(E)$ gives
    \begin{align*}
        \sum_{Q\in\mathcal{D}_\Delta(E)}\int_{\mathbb S^1} \mathcal{H}^{s-\tau}_\infty(Q_{\theta,H})\,d\mathcal{H}^1(\theta)&\lesssim \Delta^{-s}\cdot\Delta^{s-\tau}\delta^{\frac{h(1-\tau)}{2-s}-\varepsilon_0}\\
        &=\delta^{\frac{h(1-\tau)}{2-s}-k\tau-\varepsilon_0}\\
        &=\delta^{\frac{h(1-\alpha s+\alpha)}{2-s}-k\alpha(s-1)+5\left(\frac{ h}{2-s}+k\right)\varepsilon-\varepsilon_0}.
    \end{align*}
\end{proof}

\subsection{Optimizing threshold} 
In this section, we optimize the choice of $\alpha,k,h$, as introduced in the setup \ref{setup}.
\begin{prop}[Final Optimization]\label{prop:final optimizing}
    If
    \[2k+3\alpha\leqslant 1,\quad (k-\alpha)(s-1)-h\geqslant 0,\quad \text{and}\quad \frac{h(1-\alpha s+\alpha)}{2-s}-k\alpha(s-1)\geqslant 0,\]
    then
    \[\int_{\mathbb S^1} \mathcal{H}^{s-\tau}_\infty(\Vis_\theta (E))\,d\mathcal{H}^1(\theta)\lesssim_\varepsilon\delta^{\varepsilon}+\delta^{\left(\frac{h}{2-s}+k\right)\varepsilon}.\]
\end{prop}

\begin{proof}
    Recall that
    \[\Vis_\theta(E)\subseteq E_{\theta,L}\cup E_{\theta,H}\cup (\Vis_\theta(E)\cap E_{\theta,G})\cup E_{\theta,B}.\]

    For $E_{\theta,L}$, by Lemma \ref{lem:light part}, we have
    \[\mathcal{H}^{s-\tau}_\infty( E_{\theta,L})\lesssim\delta^\varepsilon,\quad\forall\,\theta\in \mathbb S^1.\]

    For $E_{\theta,H}$, we first note that, by Corollary \ref{cor:global heavy part}, we have
    \[\int_{\mathbb S^1} \mathcal{H}^{s-\tau}_\infty(E_{\theta,H,\mathscr{T}})\,d\mathcal{H}^1(\theta)\lesssim_\varepsilon\delta^\varepsilon.\]
    On the other hand, since $\frac{h(1-\alpha s+\alpha)}{2-s}-k\alpha(s-1)\geqslant 0$, by Proposition \ref{prop:local heavy part}, choosing $\left(\frac{h}{2-s}+k\right)\varepsilon\leqslant\varepsilon_0\leqslant 2\left(\frac{h}{2-s}+k\right)\varepsilon$ gives
    \[\sum_{Q\in\mathcal{D}_\Delta(E)}\int_{\mathbb S^1} \mathcal{H}^{s-\tau}_\infty(Q_{\theta,H})\,d\mathcal{H}^1(\theta)\lesssim_\varepsilon\delta^{\left(\frac{ h}{2-s}+k\right)\varepsilon}.\]
    Thus,
    \[\int_{\mathbb S^1} \mathcal{H}^{s-\tau}_\infty(E_{\theta,H})\,d\mathcal{H}^1(\theta)\lesssim_\varepsilon\delta^{\varepsilon}+\delta^{\left(\frac{h}{2-s}+k\right)\varepsilon}.\]
    
    For $\Vis_\theta(E)\cap E_{\theta,G}$, since $(k-\alpha)(s-1)-h\geqslant 0$, by Corollary \ref{cor:good part}, we have
    \[\mathcal{H}^{s-\tau}_\infty(\Vis_\theta (E)\cap E_{\theta,G})\lesssim_\varepsilon\delta^\varepsilon,\quad\forall\,\theta\in\mathbb S^1.\]
    
    For $E_{\theta,B}$, since $2k+3\alpha\leqslant 1$, by Corollary \ref{cor:bad part}, we have
    \[\int_{\mathbb S^1} \mathcal{H}^{s-\tau}_\infty(E_{\theta,B})\,d\mathcal{H}^1(\theta)\lesssim_\varepsilon\delta^\varepsilon.\]

    Combining all inequalities above gives
    \[\int_{\mathbb S^1} \mathcal{H}^{s-\tau}_\infty(\Vis_\theta (E))\,d\mathcal{H}^1(\theta)\lesssim_\varepsilon\delta^{\varepsilon}+\delta^{\left(\frac{h}{2-s}+k\right)\varepsilon}.\qedhere\]
\end{proof}

The three conditions needed are
\[2k+3\alpha\leqslant 1,\quad (k-\alpha)(s-1)-h\geqslant 0, \quad \frac{h(1-\alpha s+\alpha)}{2-s}-k\alpha(s-1)\geqslant 0.\]
The first inequality gives $k\leqslant (1-3\alpha)/2$.
The last two inequalities give
\[k\geqslant \frac{\alpha-\alpha^2(s-1)}{1-\alpha}.\]
Then we have
\[1-6\alpha+(2s+1)\alpha^2\geqslant 0.\]
In the range $0<\alpha<1/3$, the largest choice is
\[\alpha=\frac{1}{3+\sqrt{8-2s}},\quad k=\frac{1-3\alpha}{2},\quad h=\frac{(s-1)(1-5\alpha)}{2}.\]
For this choice of $\alpha$, by Proposition \ref{prop:final optimizing}, we have
\[\int_{\mathbb S^1} \mathcal{H}^{s-\tau}_\infty(\Vis_\theta (E))\,d\mathcal{H}^1(\theta)\lesssim_\varepsilon\delta^{\varepsilon}+\delta^{\left(\frac{h}{2-s}+k\right)\varepsilon}.\]
Taking $\delta\to 0$ gives \eqref{eqn:Haus content 0}.

\bibliography{refs}

@preamble{"\def\cprime{$'$} "}

@book {Mattilabook,
    AUTHOR = {Mattila, Pertti},
     TITLE = {Fourier analysis and {H}ausdorff dimension},
    SERIES = {Cambridge Studies in Advanced Mathematics},
    VOLUME = {150},
 PUBLISHER = {Cambridge University Press, Cambridge},
      YEAR = {2015},
     PAGES = {xiv+440},
      ISBN = {978-1-107-10735-9},
   MRCLASS = {28-02 (28A15 28A78 28A80 42B10 60J65)},
  MRNUMBER = {3617376},
MRREVIEWER = {Benjamin\ Steinhurst},
       DOI = {10.1017/CBO9781316227619},
       URL = {https://doi.org/10.1017/CBO9781316227619},
}

@article {Tuo23,
    AUTHOR = {Orponen, Tuomas},
     TITLE = {On the dimension of visible parts},
   JOURNAL = {J. Eur. Math. Soc. (JEMS)},
  FJOURNAL = {Journal of the European Mathematical Society (JEMS)},
    VOLUME = {25},
      YEAR = {2023},
    NUMBER = {5},
     PAGES = {1969--1983},
      ISSN = {1435-9855,1435-9863},
   MRCLASS = {28A78},
  MRNUMBER = {4592863},
       DOI = {10.4171/jems/1230},
       URL = {https://doi.org/10.4171/jems/1230},
}

@article {JJMON03,
    AUTHOR = {J\"arvenp\"a\"a, Esa and J\"arvenp\"a\"a, Maarit and MacManus,
              Paul and O'Neil, Toby C.},
     TITLE = {Visible parts and dimensions},
   JOURNAL = {Nonlinearity},
  FJOURNAL = {Nonlinearity},
    VOLUME = {16},
      YEAR = {2003},
    NUMBER = {3},
     PAGES = {803--818},
      ISSN = {0951-7715,1361-6544},
   MRCLASS = {28A80 (28A78)},
  MRNUMBER = {1975783},
MRREVIEWER = {Manuel\ Mor\'an},
       DOI = {10.1088/0951-7715/16/3/302},
       URL = {https://doi.org/10.1088/0951-7715/16/3/302},
}

@article {JJN04,
    AUTHOR = {J\"arvenp\"a\"a, Esa and J\"arvenp\"a\"a, Maarit and
              Niemel\"a, Juho},
     TITLE = {Transversal mappings between manifolds and non-trivial
              measures on visible parts},
   JOURNAL = {Real Anal. Exchange},
  FJOURNAL = {Real Analysis Exchange},
    VOLUME = {30},
      YEAR = {2004/05},
    NUMBER = {2},
     PAGES = {675--687},
      ISSN = {0147-1937,1930-1219},
   MRCLASS = {28A78 (28A15 28A80)},
  MRNUMBER = {2177426},
MRREVIEWER = {Ren\'e\ L.\ Schilling},
       URL = {http://projecteuclid.org/euclid.rae/1129416471},
}

@article {FF13,
    AUTHOR = {Falconer, Kenneth J. and Fraser, Jonathan M.},
     TITLE = {The visible part of plane self-similar sets},
   JOURNAL = {Proc. Amer. Math. Soc.},
  FJOURNAL = {Proceedings of the American Mathematical Society},
    VOLUME = {141},
      YEAR = {2013},
    NUMBER = {1},
     PAGES = {269--278},
      ISSN = {0002-9939,1088-6826},
   MRCLASS = {28A78 (28A80)},
  MRNUMBER = {2988729},
MRREVIEWER = {Ezequiel\ Rela},
       DOI = {10.1090/S0002-9939-2012-11312-7},
       URL = {https://doi.org/10.1090/S0002-9939-2012-11312-7},
}

@article {AJJRS12,
    AUTHOR = {Arhosalo, Ida and J\"arvenp\"a\"a, Esa and J\"arvenp\"a\"a,
              Maarit and Rams, Micha\l{} and Shmerkin, Pablo},
     TITLE = {Visible parts of fractal percolation},
   JOURNAL = {Proc. Edinb. Math. Soc. (2)},
  FJOURNAL = {Proceedings of the Edinburgh Mathematical Society. Series II},
    VOLUME = {55},
      YEAR = {2012},
    NUMBER = {2},
     PAGES = {311--331},
      ISSN = {0013-0915,1464-3839},
   MRCLASS = {28A78 (28A80)},
  MRNUMBER = {2928497},
MRREVIEWER = {Nicolae-Adrian\ Secelean},
       DOI = {10.1017/S0013091509001680},
       URL = {https://doi.org/10.1017/S0013091509001680},
}

@article {Dabrowski-vis24,
    AUTHOR = {D{\k{a}}browski, Damian},
     TITLE = {Visible parts and slices of {A}hlfors regular sets},
   JOURNAL = {Discrete Anal.},
  FJOURNAL = {Discrete Analysis},
      YEAR = {2024},
     PAGES = {Paper No. 17, 31},
      ISSN = {2397-3129},
   MRCLASS = {28A80 (28A78)},
  MRNUMBER = {4844555},
}

@article {DOW24-heavyslice,
    AUTHOR = {D{\k{a}}browski, Damian and Orponen, Tuomas and Wang, Hong},
     TITLE = {How much can heavy lines cover?},
   JOURNAL = {J. Lond. Math. Soc. (2)},
  FJOURNAL = {Journal of the London Mathematical Society. Second Series},
    VOLUME = {109},
      YEAR = {2024},
    NUMBER = {5},
     PAGES = {Paper No. e12910, 33},
      ISSN = {0024-6107,1469-7750},
   MRCLASS = {28A80 (28A78)},
  MRNUMBER = {4745881},
MRREVIEWER = {Wen\ Wu},
       DOI = {10.1112/jlms.12910},
       URL = {https://doi.org/10.1112/jlms.12910},
}

@article{TR26,
      title={Visibility problem in the plane}, 
      author={Tuomas Orponen and Alex Rutar},
      year={2026},
      eprint={2606.06965},
      archivePrefix={arXiv},
      primaryClass={math.CA},
      url={https://arxiv.org/abs/2606.06965}, 
    journal = {arXiv preprint, arXiv:2606.06965}
}

@article {JJSW22,
    AUTHOR = {J\"arvenp\"a\"a, Esa and J\"arvenp\"a\"a, Maarit and Suomala,
              Ville and Wu, Meng},
     TITLE = {On dimensions of visible parts of self-similar sets with
              finite rotation groups},
   JOURNAL = {Proc. Amer. Math. Soc.},
  FJOURNAL = {Proceedings of the American Mathematical Society},
    VOLUME = {150},
      YEAR = {2022},
    NUMBER = {7},
     PAGES = {2983--2995},
      ISSN = {0002-9939,1088-6826},
   MRCLASS = {28A80},
  MRNUMBER = {4428883},
MRREVIEWER = {Ying\ Xiong},
       DOI = {10.1090/proc/15843},
       URL = {https://doi.org/10.1090/proc/15843},
}

@article {ER21,
    AUTHOR = {Rossi, Eino},
     TITLE = {Visible part of dominated self-affine sets in the plane},
   JOURNAL = {Ann. Fenn. Math.},
  FJOURNAL = {Annales Fennici Mathematici},
    VOLUME = {46},
      YEAR = {2021},
    NUMBER = {2},
     PAGES = {1089--1103},
      ISSN = {2737-0690,2737-114X},
   MRCLASS = {28A80},
  MRNUMBER = {4313141},
       DOI = {10.5186/aasfm.2021.4668},
       URL = {https://doi.org/10.5186/aasfm.2021.4668},
}

@article {ONeilconnected,
    AUTHOR = {O'Neil, Toby C.},
     TITLE = {The {H}ausdorff dimension of visible sets of planar continua},
   JOURNAL = {Trans. Amer. Math. Soc.},
  FJOURNAL = {Transactions of the American Mathematical Society},
    VOLUME = {359},
      YEAR = {2007},
    NUMBER = {11},
     PAGES = {5141--5170},
      ISSN = {0002-9947,1088-6850},
   MRCLASS = {28A80 (31A15)},
  MRNUMBER = {2327025},
MRREVIEWER = {Antti\ K\"aenm\"aki},
       DOI = {10.1090/S0002-9947-07-04460-1},
       URL = {https://doi.org/10.1090/S0002-9947-07-04460-1},
}

@article{HJL-higherdim,
      title={Visibility problem in higher dimensions}, 
      author={Yuming He and Changbiao Jian and Yudong Liang},
      year={2026},
      eprint={2609.35388},
      archivePrefix={arXiv},
      primaryClass={math.CA},
    journal = {arXiv preprint, arxiv:2609.35388}
}

@article{KMP26,
      title={Visible parts and projections of the affine fractal percolation}, 
      author={Eeva Kuure and Markus Myllyoja and Aleksi Py\"or\"al\"a},
      year={2026},
      eprint={2610.21768},
      archivePrefix={arXiv},
      primaryClass={math.CA},
    journal = {arXiv preprint, arxiv:2610.21768}
}
\bibliographystyle{abbrv}

\end{document}